\documentclass{amsart}
\usepackage{amssymb, amsmath, amsfonts, amscd}

\input xypic

\theoremstyle{plain}
\newtheorem{theorem}{Theorem}[section]
\newtheorem{corollary}[theorem]{Corollary}
\newtheorem{lemma}[theorem]{Lemma}
\newtheorem{proposition}[theorem]{Proposition}

\newtheorem{example}[theorem]{Example}

\theoremstyle{definition}
\newtheorem{definition}[theorem]{Definition}

\theoremstyle{remark}

\numberwithin{equation}{theorem}

\newcommand{\Hom}{\operatorname{Hom} }

\newcommand{\D}{\Delta}
\newcommand{\tD}{\tilde{\Delta}}

\newcommand{\DV}{\Delta_V }
\newcommand{\DR}{\Delta_R }
\newcommand{\DVd}{\Delta_{V^*} }

\newcommand{\Sym}{\operatorname{S} }

\newcommand{\e}{\epsilon}

\newcommand{\co}{\underline{comod}}

\newcommand{\Spec}{\operatorname{Spec} }

\newcommand{\SLA}{\operatorname{SL}_{2,A} }
\newcommand{\SLZ}{\operatorname{SL}_{2,\mathbb{Z}}}
\newcommand{\SL}{\operatorname{SL}}

\newcommand{\Dist}{\operatorname{Dist} }
\newcommand{\GL}{\operatorname{GL} }
\newcommand{\Lie}{\operatorname{Lie} }

\newcommand{\alg}{\underline{A-alg}}

\newcommand{\zalg}{\underline{\mathbb{Z}-alg}}

\newcommand{\grp}{\underline{Groups}}
\newcommand{\Z}{\mathbb{Z}} 
\newcommand{\T}{\tilde{T}}

\begin{document}

\title{A universal Clebsch-Gordan filtration for   $\GL_{2,A}$ }

\author{Helge \"{O}ystein Maakestad }





\date{ august 2023} 

\begin{abstract}  In this paper I prove a version of the Clebsch-Gordan formula in the following sense: Let $\Z$ be the ring of integers and let $V$ be the standard $\SL_{2,\Z}$-comodule.
There is for any pair of integers $n,m$ with $1 \leq n \leq m$ a finite filtration $ F_i \subseteq \Sym^n(V)\otimes \Sym^m(V)$ of comodules on $\SL_{2,\Z}$  with $F_i/F_{i+1} \cong \Sym^{n+m -2i}(V)$. Hence in the grothendieck group of finite rank comodules we get an equality of classes - a virtual Clebsch-Gordan formula:
\[  [\Sym^n(V)] [ \Sym^m(V)] = \sum_{i=0}^n [\Sym^{n+m-2i}(V)]. \] 

The filtration  $F_i$ pulls back to a filtration $F_i \otimes_{\Z} A$ of $\Sym^n(V_A) \otimes \Sym^m(V_A)$ for any commutative ring $A$, and it gives a virtual Clebsch-Gordan formula in the Grothendieck group of "finite rank" comodules on $\SLA$ for any such $A$. I prove a similar result for $\GL_{2,A}$  giving an explicit construction of a finite filtration of $\GL_{2,A}$-comodules   $E_i \subseteq \Sym^n(V_A) \otimes \Sym^m(V_A)$  for any commutative ring $A$ with similar properties: There is an isomorphism of comodules $ E_i/E_{i+1} \cong \Sym^{n+m-2i}(V_A) \otimes (\wedge^2 V_A)^{\otimes_A^i}$ and a virtual Clebsch-Gordan formula

\[  [\Sym^n(V_A)] [ \Sym^m(V_A)] = \sum_{i=0}^n [\Sym^{n+m-2i}(V_A)][\wedge^2 V_A)]^i. \] 

 I introduce the notion "good filtration" for a torsion free comodule on $\SL_{2,A}$ and $\GL_{2,A}$ and give non trivial examples of torsion free modules equipped with a good filtration. The filtrations $F_i, E_i$ are examples of good filtrations. 
I give one "constructive proof" proof using the Hopf algebra $k[\GL_{2,A}]$ and "local coordinates"  and one "non-constructive proof" proof using the functor of points $h_{\GL_{2,A}}(-)$.
\end{abstract}

\maketitle
\tableofcontents

\section{Introduction}  The aim of the paper is to study torsion free finite rank comodules on the group scheme $\SLA$ where $A$ is a Dedekind domain. If $A$ is a field of characteristic zero there is a complete classification of finite dimensional irreducbible comodules. For any such comodule $W$ there is a direct sum decomposition $W \cong \oplus_{i=1}^n \Sym^{d_i}(V)$ where $V$ is the "standard comodule" and $\Sym^d(V)$ is the d'th symmetric power of $V$. Moreover for any integer $d \geq 0$ it follows the comodule $\Sym^d(V)$ is an irreducible comodule. These properties no longer hold for a general Dedekind domain such as $\Z, \Z_{(p)}$ or $\Z[\frac{1}{m}]$. For a torsion free finite rank comodule $W$ on $\SLA$ 
there is no direct sum decomposition into symmetric powers of the standard comodule $V$ (see Theorem \ref{completereducibility}). Moreover the symmetric powers $\Sym^d(V)$ are no longer irreducible comodules on $\SLA$ (see Example \ref{notirreducible}). 

We say a torsion free comodule $W$ has a "good filtration"     
of comodules iff there is a finite filtration of comodules $0:=F_{n+1} \subsetneq F_n \subsetneq \cdots \subsetneq F_0:=W$ 
with the following property: There is an isomorphism of comodules $F_i/F_{i+1} \cong \Sym^{d_i}(V)$ with $d_i \geq 0$ an integer for each $i$.  When $A$ is a field of chcaracteristic zero there is the "classical Clebsch-Gordan formula". It is the following formula:   Let $1 \leq n \leq m$ be integers. There is an isomorphism

\[  \Sym^n(V) \otimes \Sym^m(V) \cong  \oplus_{i=0}^n \Sym^{n+m-2i}(V)  \]

of comodules on $\SLA$. This formula gives many filtrations  of the tensor product. One of them is the following filtration: Let  $F(k)_i := \oplus_{j=0}^n \Sym^{n+m-2j}(V)$ It follows we get a filtration -  \emph{the classical Clebsch-Gordan filtration} -

\[   0:=F(k)_{n+1} \subsetneq F(k)_n \subseteq F(k)_1 \subseteq F(k)_0:= \Sym^n(V) \otimes \Sym^m(V) \]

of modules. The classical Clebsch-Gordan filtration is an example of a good filtration of $\Sym^n(V)\otimes \Sym^m(V)$.
In this paper we study the tensor product $\Sym^n(V)\otimes \Sym^m(V)$ when $A$ is a Dedekind domain. There is a "canonical"  filtration
called \emph{the universal Clebsch-Gordan filtration}:

\[  0:=F_{n+1} \subsetneq  F_n \subseteq F_{n-1} \subseteq \cdots \subseteq F_1 \subseteq F_0:= \Sym^n(V) \otimes \Sym^m(V) .\]

This filtration is a filtration of $\SL_{2,A}$-comodules with $F_i/F_{i+1} \cong \Sym^{n+m-2i}(V)$ for all $i=0,..,n$ (see Corollary \ref{CGfiltration}). The filtration $F_i$ pulls back to the "classical Clebsch-Gordan filtration"  filtration in the case when $k$ is a field of characteristic zero. 
It follows the universal Clebsch-Gordan filtration is a good filtration of comodules of $\Sym^n(V) \otimes \Sym^m(V)$ for the group scheme $\SLA$ over any commutative ring $A$. Hence the Clebsch-Gordan decomposition does not generalize
to an arbitrary commutative ring but the classical Clebsch-Gordan filtration does.

Similarly I prove there is a filtration 

\[0:= E_{n+1} \subsetneq E_n \subsetneq \cdots \subsetneq E_1 \subseteq E_0:=\Sym^n_A(V) \otimes \Sym^m_A(V) \]

called the \emph{universal Clebsch-Gordan filtration} for the $\GL_{2,A}$-comodule $\Sym_A^n(V)\otimes \Sym_A^m(V)$.

There is for every $i$ a canonical isomorphism of comodules

\[ E_i/E_{i+1} \cong \Sym_A^{n+m-2i}(V) \otimes (\wedge^2_A V)^{\otimes_A i} .\]

In the Grothendieck group of finite rank $\GL_{2,A}$-comodules we get an equality of classes

\[   [\Sym^n_A(V)][\Sym_A^m(V)] =  \sum_{i=0}^n [\Sym_A^{n+m-2i}(V) ][\wedge^2_A V]^i ,\]

called the \emph{virtual Clebsch-Gordan formula} for $\GL_{2,A}$.

We study the "complete reducibility property" for group schemes $\SLA$ where $A:=\Z, \Z_{(p)}, \Z[\frac{1}{m}]$, and find that these group schemes are not "completely reducible" in the sense that there are surjections of comodules
 $\pi: V \rightarrow W$ that do not have a section (see Theorem \ref{completereducibility}).  Moreover: The comodule $\Sym^n(V)$ is no longer irreducible when $A$ is one of the rings $\Z, \Z_{(p)}$ or $ \Z[\frac{1}{m}]$. Hence there is no classification of irreducible $\SLA$-comodules in terms of direct sums of the symmetric powers $\Sym^n(V)$ of the standard module $V$.

We give explicit examples of comodules $W, W^*$ on $\SLZ$ with the property that $W$ and $W^*$ have the same weight space decomposition but where $W$ and $W^*$ are non isomorphic as comodules. There is an isomorphism
$W \otimes_{\Z} \Z[\frac{1}{2}] \cong W^* \otimes_{\Z} \Z[\frac{1}{2}]$ as comodules on $\SL_{2,\Z[\frac{1}{2}]}$. Hence the modules $W, W^*$ are not classified by their weights and weight spaces (see Theorem \ref{descent}).
We also give examples of comodules with no good filtration (see Example \ref{nogoodfiltration}).

Note: In a paper from 2020 (see \cite{maa1}) I initiated a study of flag bundles $\mathbb{F}(E)$ where $E$ is a finite rank vector bundle on $S:=\Spec(A)$ and where $A$ is a Dedekind domain. If $P \subseteq \SL_{n,A}$ is a parabolic subgroup scheme of 
$\SL_{n,A}$, it follows the quotient $\SL_{n,A}/P$ is such a flag bundle. Any representation $(V, \rho)$ of $P$ with the property that $V$ is a torsion free and finite rank projective $A$-module, gives rise to a finite rank vector bundle $E(\rho)$ on $\SL_{n,A}/P$
and for this reason one wants to be able to classify such representations. It is hoped the constructions introduced in this paper may have applications in the study of such flag bundles and their vector bundles.

\section{Symmetric powers and the universal Clebsch-Gordan filtration}

In this section we wil study representations of the group scheme $\SLA$ where $A$ is any commutative ring and $\SLZ$ where $\Z$ is the ring of integers. We will us the language of "affine group schemes" introduced in the books of Demazure, Gabriel and Jantzen (see \cite{DG} and \cite{jantzen}). If $G:=\Spec(R)$ is an affine group scheme with Hopf algebra $R$ we will let $h_R$ denote the "functor of points" of the group scheme $G$.

Let in the following $A$ be a fixed commutative unital ring and let $V:=A\{e_1,e_2\}$ be the free $A$-module on the elements $e_1,e_2$. Let $x_{11}, x_{12}, x_{21}, x_{22}$ be independent variables over $A$ and let $D:=x_{11}x_{22}-x_{12}x_{21}$
Let $R:=A[x_{ij}]/(D-1)$ and let $\SLA:=\Spec(R)$. It follows the functor of points 

\[  h_R: \alg \rightarrow \grp \]

defined by $h_R(B):=\Hom_{A-alg}(R,B)$ where $B$ is an $A$-algebra, has the property that $h_R(B)$ equals the group of $2 \times 2$-matrices with coefficients in $B$ and determinant $1$. Hence the functor $h_R$ is an affine group scheme in the sense of Jantzen's book \cite{jantzen} on algebraic groups. Since $R$ gives rise to an affine group scheme it follows the ring $R$ has the structure as Hopf algebra $(R, \Delta_R, S, \epsilon)$. Here $\Delta_R, S, \epsilon$ are maps

\[  \Delta_R: R \rightarrow R \otimes_A R ,\]

\[ S: R \rightarrow R\]

and

\[ \epsilon: R \rightarrow A \]

satisfying a set of criteria. The map $\Delta_R$ is the comultiplication, $S$ is the co-inversion and $\epsilon$ is the counit. The map $\Delta_R$ is defined as follows

\[ \Delta_R(x_{11}):= x_{11} \otimes x_{11} + x_{12} \otimes x_{21}, \Delta_R(x_{12}):=x_{11} \otimes x_{12} + x_{12} \otimes x_{22} ,\]

\[ \Delta_R(x_{21}):= x_{21} \otimes x_{11} + x_{22} \otimes x_{21}, \Delta_R(x_{22}):= x_{21} \otimes x_{12} + x_{22} \otimes x_{22}.\]

Moreover 

\[ S(x_{11}):= x_{22}, S(x_{12}):=-x_{12}, S(x_{21}):=-x_{21}, S(x_{22}):=x_{11},\]

and

\[ \epsilon(x_{11}):=\epsilon(x_{22}):=1, \epsilon(x_{12}):=\epsilon(x_{21}):=0.\]

One checks the above formulas define a Hopf algebra in the sense of Jantzen's book. Define the following comodule structure on $V$:

\[ \Delta_V: V \rightarrow R \otimes V \]

by

\[  \DV(e_1):= x_{11} \otimes e_1 + x_{12} \otimes e_2, \DV(e_2):= x_{21} \otimes e_1 + x_{22} \otimes e_2.\]

One checks there are equalitites

\[ 1\otimes \DV \circ \DV = \Delta_R \otimes 1 \circ \DV \]

and 

\[ \epsilon \otimes 1 \circ \DV = I_V.\]

Hence $(V, \DV)$ is a comodule on $\SLA$.

\begin{definition} \label{standard} Let $(V, \DV)$ be the \emph{standard comodule} on $\SLA$.
\end{definition}

\begin{example} A relation to multiplication with a matrix. \end{example}

The structure $\DV$ as comodule on $V$ corresponds to "matrix multiplication from the right". Let $v\in V$ be represented by a pair $[a_1,a_2]$ with $a_i \in A$,
 and choose a $A$-valued point $M \in \SLA(A)$, where $M$  is a $2 \times 2$-matrix with coefficients in $A$
and determinant one. The matrix $M$ acts on the vector $[a_1,a_2]$ using matrix multiplication from the right.
There is another obvious  comodule structure

\[ \D^l_V: V \rightarrow V \otimes R \]

where we use matrix multiplication from the left.

\begin{example} \label{naivedual} Symmetric powers, exterior powers and "classical duals" \end{example}

We may for any integer $n \geq 0$ construct the symmetric power $\Sym^n(V) \cong \oplus_{i=0}^n A e_1^{n-i}e_2^i$. It follows there is a canonical structure as comodule

\[  \Delta_n: \Sym^n(V) \rightarrow R \otimes \Sym^n(V) .\]

We may take the exterior power $\wedge^2 V$ and there is a canonical comodule structure

\[ \Delta_1: \wedge^2 V \rightarrow R \otimes \wedge^2 V. \]

One checks this is the "trivial" structure. Let $v:=e_1 \wedge e_2$ be a basis for the free rank one $A$-module $\wedge^2 V$. It follows $\Delta_1(av)=1 \otimes av$ for any  $a \in A$. 

Let the "classical dual" $V^*:=A\{x_1,x_2\}$ be the free $A$-module on the basis $x_i:=e_i^*$ and define the following map:

\[ \DVd: V^* \rightarrow R \otimes V^* \]

by

\[  \DVd(x_1):= x_{22} \otimes x_1 - x_{21} \otimes x_2, \DVd(x_2):= -x_{12} \otimes x_1 + x_{11} \otimes x_2.\]

One checks $(V^*, \DVd)$ is a comodule. 

\begin{lemma} There is an isomorphism of comodules $\phi: V \rightarrow V^*$ defined by $\phi(e_1):= x_2, \phi(e_2):=-x_1$.
\end{lemma}
\begin{proof} One checks the map $\phi$ is a map of comodules which is an isomorphism. The Lemma follows.
\end{proof}

Hence for any integer $n \geq 0$ it follows there is an isomorphism $\Sym^n(V) \cong \Sym^n(V^*)$ of comodules. Over a field $k$ this can be deduced from the fact that the canonical map

\[  f: V \otimes_k V \rightarrow \wedge^2 V \]

defined by $f(u \otimes v):= u \wedge v$ is a nondegenerate pairing and the fact that $\wedge^2 V \cong T$ is the trivial comodule. This result is by the above argument true for $\SLA$ over any commutative ring $A$.

Let in the following $A:=\Z$ be the ring of integers. We let $\Sym^*(V):= \oplus_{d \geq 0} \Sym^d(V) \cong \Z[e_1,e_2]$ be the symmetric algebra of $V$ over $\Z$. This is a graded commutative ring with graded pieces $\Sym^d(V)$ - homogeneous polynomials in $e_1,e_2$ of degree $d$.

Since $G:=\SLZ$ is a group scheme over $\Z$ it follows its algebra of distributions $\Dist(G)$ acts on $V$ and $\Sym^d(V)$. There is an inclusion $\Lie(G) \subseteq \Dist(G)$ and $\Lie(G) \cong \mathfrak{sl}(2,\Z)$
is the Lie algebra of $2\times 2$-matrices with coefficients in $\Z$ and trace equal to zero. It follows $\Sym^d(V)$ is a module on $\Lie(G)$ and it has by \cite{DG}  a weight space decomposition for the action of $H$. We let $\Lie(G):=\Z\{x,H,y\}$
be the "standard generators" for $\Lie(G)$.   The elements $x,H,y$  are defined as follows:

\[
x:=  
\begin{pmatrix}   0  &  1  \\
           0  &  0 
\end{pmatrix}.
\]
\[
H:=  
\begin{pmatrix}   1  &  0  \\
           0  &  -1 
\end{pmatrix}.
\]

and

\[
y:=  
\begin{pmatrix}   0  &  0  \\
           1  &  0 
\end{pmatrix}.
\]

If $v:=e_1^pe_2^q \in \Sym^d(V)$ it follows

\[  H(v)=(p-q) v.\]

If we let $V(d):= \Sym^d(V)$ it follows there is by \cite{DG} a direct sum decomposition into eigenspaces for the action of $H$:

\[  V(d) \cong \oplus_{i=0}^d V(d)_{-d+2i}  \]

where $ V(d)_j$ is the eigenspace of vectors where $H$ acts by multiplying with $j$. As an example: $H(e_2^d)= -d e_2^d$. Hence the "vector" $v:=e_2^d$ has weight equal to $-d$.

There is a sequence of maps  of $\Z$-modules

\[   \Sym^{n-1}(V) \otimes \Sym^{m-1}(V) \rightarrow^{f_{n,m}} \Sym^n(V) \otimes \Sym^m(V) \rightarrow^{g_{n,m}} \Sym^{n+m}(V) \]

where $g_{n,m}$ is the multiplication map and $f_{n,m}$ is defined as follows:

\[ f_{n,m}(e_1^pe_2^q \otimes e_1^se_2^t):=  e_1^{p+1}e_2^q \otimes e_1^se_2^{t+1}  -  e_1^pe_2^{q+1} \otimes e_1^{s+1}e_2^t.\]

For $W:= \Sym^{n+m}(V)$ there is a weight space decomposition $W \cong \oplus_{i=0}^{n+m} W_{-(n+m) +2i}$  and the weight of $v:=e_1^pe_2^q$ denoted $wt(v)$ equals $p-q$. It follows

$wt(e_1^ie_2^{n+m-i})=2i-(n+m)$. Let $0 \leq p \leq n, 0 \leq s \leq m$ and let

\[ v_k:= e_1^ke_2^{n-k} \otimes e_1^{i-k}e_2^{m-i+k} .\]

It follows

\[ f_{n,m}(v_k)=e_1^ie_2^{m+n-i}.\]

Let $v:= e_1^pe_2^{n-p} \otimes e_1^se_2^{m-s}$. It follows

\[ f_{n,m}(v)=e_1^{p+s}e_1^{n+m-(p+s)}.\]

It follows $wt(f_{n,m}(v))=wt(f_{n,m}(v_k))$ iff $ p+s =i$. It follows the pair $(p,s)$ looks as follows:

\[ (p,s)=\{(i,0), (i-1,1),\ldots, (k,i-k),\ldots, (0,i)\}.\]

Since $(q,t)=(n-p, m-s)$ it follows the corresponding 4-tuple $(p,q,s,t)$  must look as follows:

\[ ( p,q, s,t)=(k, n-k, i-k, m-i+k)\]

for an integer $0 \leq k \leq i$. Hence  a vector $v:= e_1^pe_2^{n-p} \otimes e_1^se_2^{m-s}$ has the property that $f_{n,m}(v)=e_1^ie_2^{n+m-i}$  iff there is an integer $k$ with  $0 \leq k \leq i$ and with 

\[ v_k:=e_1^ke_2^{n-k} \otimes e_1^{i-k}e_2^{m-i+k}.\]

Hence an element $v:= \sum_{j=0}^i a_iv_i$ has $f(v)=0$ iff $a_0+a_1+\cdots +a_i=0$. For each $1 \leq k \leq i-1$ it follows

\[ v_k -v_{k+1} = -z e_1^ke_2^{n-(k+1)} \otimes e_1^{i-(k+1)}e_2^{m-(i-k)} \]

where $z:=e_1\otimes e_2 - e_2 \otimes e_1$.

\begin{lemma} \label{zmodules}  There is an exact sequence of free $\Z$-modules

\[ 0 \rightarrow \Sym^{n-1}(V) \otimes \Sym^{m-1}(V) \rightarrow^{f_{n,m}} \Sym^n(V) \otimes \Sym^m(V)  \rightarrow^{g_{n,m}} \Sym^{n+m}(V) \rightarrow 0\]

where $g_{n,m}$ is the canonical "multiplication map" and $f_{n,m}$ is defined as follows: Let $e_1^pe_2^q \otimes e_1^se_2^t$ be an element with $p+q=n-1, s+t=m-1$. Define $f_{n,m}$ as follows:

\[ f_{n,m}(e_1^pe_2^q \otimes e_1^se_2^t):=  e_1^{p+1}e_2^q \otimes e_1^se_2^{t+1}  -  e_1^pe_2^{q+1} \otimes e_1^{s+1}e_2^t.\]

\end{lemma}
\begin{proof} It is clear that  $Im(f_{n,m}) \subseteq ker(g_{n,m})$. We prove the reverse inclusion. By the above calculation it follows a vector $v:=a_0v_0+\cdots a_i v_i $ has the property that $f_{n,m}(v)=0$ iff
$a_0 + \cdots + a_i=0$. We may write

\[ v =a_0(v_0-v_1)+(a_0+a_1)(v_1-v_2) + \cdots (a_0+ \cdots a_{i-1})(v_{i-1}-v_i)+ (a_0+\cdots + a_i)v_i =\]

\[ a_0(v_0-v_1)+ \cdots + (a_0+\cdots + a_{i-1})(v_{i-1}-v_i) \]

and it follows $v \in Im(f_{n,m})$ since each vector $v_k-v_{k+1}$ is in $Im(f_{n,m})$. The Lemma follows.

\end{proof}

Note: The result in Lemma \ref{zmodules} holds if we replace $\Z$ by any commutative unital ring $A$.

 Let $\Sym^*(V):= \oplus_{n \geq 0} \Sym^n(V)$ be the symmetric algebra of $V$ over $\Z$.  Let $W:= \Sym^*(V) \otimes \Sym^*(V):=\oplus_{d \geq 0} W_d$ be the graded ring with grading

\[ W_d:= \oplus_{i+j=d} \Sym^i(V) \otimes \Sym^j(V) \]

It follows $\Sym^*(V) \cong \Z[e_1,e_2]$ is the polynomial ring on $e_1,e_2$ and $W:=\Z[e_1,e_2]\otimes_{\Z} \Z[e_1,e_2] \cong \Z[e_1,e_2,u_1,u_2]$ is the polynomial ring in 4 independent variables.

\begin{lemma} There are comodule structures

\[ \D_1: \Sym^*(V) \rightarrow R \otimes \Sym^*(V) \]

and

\[ \D_2:  W \rightarrow R \otimes W \]

and the canonical multiplication map $g: W \rightarrow \Sym^*(V)$ is a map of comodules.
\end{lemma}
\begin{proof} Define the map $\D_1$ as follows: Let $v:=e_1^p e_2^q \in \Sym^n(V)$ be two elements and define

\[ \D_1(v):= \DV(e_1)^p\DV(e_2)^q  \]

and for any polynomial $f(e_1,e_2) \in \Sym^*(V)$ define

\[  \D_1(f):= f(\DV(e_1), \DV(e_2)) .\]

By definition $\Delta_1, 1 \otimes \Delta_1$  and $\DR \otimes 1$ are maps of $\Z$-algebras. It follows 

\[ 1\otimes \Delta_1(\Delta_1(v))= \]

\[  1\otimes \Delta_1(\DV(e_1)^p\DV(e_2)^q))= \]

\[  ((1 \otimes \Delta_1(x_{11} \otimes e_1+x_{12} \otimes e_2))^p(1\otimes \Delta_1(x_{21} \otimes e_1 + x_{22} \otimes e_2))^q=\]

\[ (x_{11} \otimes \DV(e_1) + x_{12}\otimes \DV(e_2))^p(x_{21} \otimes \DV(e_1) + x_{22} \otimes \DV(e_2))^q =\]

\[  (\DR(x_{11}) \otimes e_1 + \DR(x_{12}) \otimes e_2)^p(\DR(x_{21} \otimes e_1 + \DR(x_{22}) \otimes e_2)^q =\]

\[ (\DR \otimes 1(\DV(e_1)))^p(\DR \otimes 1(\DV(e_2)))^q =\]

\[ \DR \otimes 1 (\Delta_1(e_1)^p \Delta_1(e_2)^q) = \DR \otimes 1(\Delta_1(e_1^p e_2^q)) .\]

It follows using a similar calculation that $\D_1: \Sym^*(V) \rightarrow R \otimes \Sym^*(V)$ is a map of $\Z$-modules satisfying the axioms of being a comodule and the first claim follows.  Note: The map $\Delta_1$ 
is a map of $\Z$-algebras.








 In a similar way one checks the tensor product $W:= \Sym^*(V) \otimes_{\Z} \Sym^*(V)$ has a comodule structure defined as follows: Let $v \otimes w \in W$ be an element with 

\[  \Delta_1(v) := \sum_i x_i \otimes u_i, \Delta_1(w):= \sum_j y_j \otimes v_j.\]

The differential $\Delta_2: W \rightarrow R \otimes W$ is defined by $\Delta_2(v \otimes w):= \sum_{i,j} x_iy_j \otimes u_i \otimes v_j$. Define a map

$g: W \rightarrow \Sym^*(V)$ by $g(v \otimes w):= vw$. It follows 

\[  1 \otimes g(\Delta_2(v \otimes w))= \sum_{i,j} x_iy_j \otimes u_iv_j \]

and

\[ \Delta_1(g(v\otimes w))= \Delta_1(vw):=\Delta_1(v)\Delta_1(w)= \sum_{i,j} x_iyj \otimes u_iv_j \]

since $\Delta_1$ is a map of rings. It follows $g$ is a map of comodules and the Lemma follows.
\end{proof}

By definition the map $\D_1$ is defined as follows: Let $f(e_1,e_2) \in \Sym^*(V)$ be a polynomial. We get

\[ \D_1(f):= f(\DV(e_1),\DV(e_2)). \]  

Let $W:=\Sym^*(V) \otimes \Sym^*(V)$ and let $v \otimes w \in W$ with $\D_1(v)=\sum_i x_i \otimes v_i, \D_1(w):= \sum_j y_j \otimes u_j$. It follows by definition

\[  \D_2(v\otimes w):= \sum_{i,j} x_iy_j \otimes v_i \otimes u_j \in R \otimes W.\]

Let $z:= e_1 \otimes e_2 - e_2 \otimes e_1 \in W$.  By definition we have

\[ \DV(e_1):= x_{11} \otimes e_1 + X_{12} \otimes e_2, \DV(e_2):= x_{21} \otimes e_1 + x_{22} \otimes e_2.\]  It follows

\[ \phi_z(v\otimes w):= z(v\otimes w)= e_1 v \otimes e_2 w- e_2 v \otimes e_1 w. \]

The map $\D_2$ is defined as follows:  There is a canonical map

\[ \eta: R \otimes \Sym^*(V) \otimes R \otimes \Sym^*(V) \rightarrow R \otimes W \]

defined by 

\[ \eta( x \otimes u  \otimes y \otimes v):= xy \otimes u \otimes v.\]

Let $\rho:= \D_1 \otimes \D_1$.
By definition it follows $\D_2:= \eta \circ \rho$. We want to calculate $\D_2(\phi_z(v \otimes w))$.  We get

\[  \rho(\phi_z(v\otimes w)):= \rho( e_1 v \otimes e_2 w- e_2 v \otimes e_1 w)= \]

\[  \D_1(e_1v)\otimes\D_1(e_2w)- \D_1(e_2v)\otimes \D_1(e_1w) =\]

\[ \DV(e_1)\D_1(v)\otimes \DV(e_2)\D_1(w)- \DV(e_2)\D_1(v) \otimes \DV(e_1)\D_1(w).\]

It follows 

\[ \rho(e_1v \otimes e_2w)=  \D_1(e_1v)\otimes\D_1(e_2w)= \DV(e_1)\D_1(v)\otimes \DV(e_2)\D_1(w) =\]

\[  \sum_{i,j} x_{11}x_i \otimes e_1u_i \otimes x_{21}y_j \otimes e_1v_j + x_{11}x_i \otimes e_1u_i \otimes x_{22}y_j \otimes e_2v_j +\]

\[  x_{12}x_i \otimes e_2u_i \otimes x_{21}y_j \otimes e_1v_j + x_{12}x_i \otimes e_2u_i \otimes x_{22}y_j \otimes e_2v_j.\]

It follows

\[ \eta(\rho( e_1v \otimes e_2 w))= \]

\[  \sum_{i,j} x_{11}x_{21}x_iy_j \otimes e_1u_i \otimes e_1v_j + x_{11}x_{22}x_iy_j  \otimes e_1u_i  \otimes e_2v_j +\]

\[  x_{12}x_{21} x_iy_j  \otimes e_2u_i \otimes e_1v_j + x_{12}x_{22}x_iy_j  \otimes e_2u_i  \otimes e_2v_j.\]

Similarly it follows 

\[ \eta(\rho(e_2v \otimes e_1 w))= \]

\[ \sum_{i,j} x_{11}x_{21}x_iy_j \otimes e_1u_i \otimes e_1v_j + x_{12}x_{21}x_iy_j  \otimes e_1u_i  \otimes e_2v_j +\]

\[  x_{11}x_{21} x_iy_j  \otimes e_2u_i \otimes e_1v_j + x_{12}x_{22}x_iy_j  \otimes e_2u_i  \otimes e_2v_j.\]

It follows

\[ \D_2(\phi_z(v\otimes w))= \]

\[ \sum_{i,j} x_i y_j \otimes ((e_1 \otimes e_2 - e_2 \otimes e_1)u_i \otimes v_j) =\]

\[ 1\otimes z(\sum_{i,j} x_iy_j \otimes u_i \otimes v_j)=\]

\[ 1 \otimes \phi_z( \sum_{i,j} x_iy_j \otimes u_i \otimes v_j)= 1 \otimes \phi_z (\D_2(v \otimes w)).\]

\begin{proposition} \label{comodules} Let $W:=\Sym^*(V) \otimes \Sym^*(V)$ and let $z:= e_1 \otimes e_2 -e_2 \otimes e_1 \in W$. Define the map $\phi_z: W\rightarrow W$  by $\phi_z(u):= zu$. It follows $\phi_z$ is a map of comodules.

\end{proposition}
\begin{proof} Let $v \otimes w \in W$ with $\D_1(v)=\sum_i x_i \otimes v_i, \D_1(w):= \sum_j y_j \otimes u_j$. It follows by definition

\[  \D_2(v\otimes w):= \sum_{i,j} x_iy_j \otimes v_i \otimes u_j \in R \otimes W.\]

Let $z:= e_1 \otimes e_2 - e_2 \otimes e_1 \in W$.  By the above calculation it follows $\D_2(\phi_z(v\otimes w))= 1\otimes z(\D_2(v\otimes w))=1\otimes \phi_z(\D_2(v\otimes w))$. The Proposition follows.

\end{proof}

\begin{theorem} \label{exactthm} Let $1 \leq n \leq m$ be a pair of integers. There is an  exact sequence of $\SLZ$-comodules

\[ 0 \rightarrow \Sym^{n-1}(V) \otimes \Sym^{m-1}(V) \rightarrow^{f_{n,m}} \Sym^n(V) \otimes \Sym^m(V)  \rightarrow^{g_{n,m}} \Sym^{n+m}(V) \rightarrow 0.\]

\end{theorem}
\begin{proof} The proof follows from Lemma \ref{zmodules} and Proposition \ref{comodules}. If $1\leq n \leq m$ are integers  it follows the map $\phi_z$ induce
 a map $\phi^{n,m}_z$ at the graded pieces:

\[ \phi^{n,m}_z: W_{n+m-2} \rightarrow W_n .\]

For the module $\Sym^{n-1}(V) \otimes \Sym^{m-1}(V)$ we get the map

\[ f_{n,m}:  \Sym^{n-1}(V) \otimes \Sym^{m-1}(V) \rightarrow \Sym^{n}(V) \otimes \Sym^{m}(V) \] 

from Lemma \ref{zmodules}. The Theorem follows.
\end{proof}

\begin{corollary}\label{CGfiltration}  Let $1 \leq n \leq m$ be a pair of integers. There is a filtration of $\SLZ$-comodules 

\[  0:=F_{n+1} \subsetneq  F_n \subseteq F_{n-1} \subseteq \cdots \subseteq F_1 \subseteq F_0:= \Sym^n(V) \otimes \Sym^m(V) \]

with isomorphisms $F_i/F_{i+1} \cong \Sym^{n+m-2i}(V)$  for $i=0,..,n$.

\end{corollary}
\begin{proof} The Corollary follows from Theorem \ref{exactthm}.

\end{proof}

\begin{definition} Let $1 \leq n \leq m$ be a pair of integers. The filtration of $\SLZ$-comodules 

\[  0:=F_{n+1} \subsetneq  F_n \subseteq F_{n-1} \subseteq \cdots \subseteq F_1 \subseteq F_0:= \Sym^n(V) \otimes \Sym^m(V) \]

from Corollary \ref{CGfiltration} is the \emph{universal Clebsch-Gordan filtration}.
\end{definition}

\begin{corollary} \label{virtualCG} (The virtual  Clebsch-Gordan formula) Let $1 \leq n \leq m$ be a pair of integers. Let  $\operatorname{K}(\SLZ)$ be the grothendieck group of finite rank comodules on $\SLZ$.
There is an equality

\[  [\Sym^n(V)] [ \Sym^m(V)] = \sum_{i=0}^n [\Sym^{n+m-2i}(V)]. \] 

in $\operatorname{K}(\SLZ)$. 
\end{corollary} 
\begin{proof}  The Corollary follows from Theorem \ref{exactthm}.
\end{proof}

\begin{example}  The classical Clebsch-Gordan formula. \end{example}

Let $S.=\Spec(\Z)$ and let $\pi: \SLZ:=G \rightarrow S$ be the canonical map. Let $k$ be the field of rational numbers and let $T:=\Spec(k)$. It follows there is an isomorphism $G_k:=\pi^{-1}(T) \cong \SL(2,k)$ and we may view
$G_k$ as the generic fiber of the map $\pi$. The group scheme $G_k$ is completely reducible and we may pull back $\Sym^n(V) \otimes \Sym^m(V)$ to $G_k$ to get the comodule $\Sym^n(V_k) \otimes_k \Sym^m(V_k)$ where $V_k:=V \otimes_{\Z} k$.
The Clebsch-Gordan formula holds for $\Sym^n(V_k) \otimes \Sym^m(V_k)$:  Let $1 \leq n \leq m$ be a pair of integers. There is an isomorphism of $G_k$-comodules

\[  \Sym^n(V_k) \otimes \Sym^m(V_k) \cong \oplus_{j=0}^n \Sym^{n+m-2j}(V_k).\]

We get a filtration  of $G_k$-comodules - the \emph{classical Clebsch-Gordan filtration} -

\[  F(k)_i:= \sum_{j=i}^n \Sym^{n+m-2j}(V_k) \]

with  $F(k)_i/F(k)_{i+1} \cong \Sym^{n+m-2i}(V_k)$.  The filtration of comodules

\[   0:=F(k)_{n+1} \subsetneq F(k)_n \subseteq F(k)_1 \subseteq F(k)_0:= \Sym^n(V_k) \otimes \Sym^m(V_k) \]

is the pull back of the filtration $F_i \subseteq \Sym^n(V) \otimes \Sym^m(V)$ from Lemma \ref{CGfiltration} in the sense there are isomorphisms of comodules

\[   \pi^*(F_i) \cong F(k)_i \]

for all $i=0,..,n$. Hence the Clebsch-Gordan formula does not generalize to $\SLZ$ but the filtration does.  There are many different filtrations of the comodule $\Sym^n(V_k) \otimes \Sym^m(V_k)$. The filtration $F(k)_i$ gives an ordering
of the direct summands $\Sym^{n+m-2i}(V_k)$ of $\Sym^n(V_k)\otimes \Sym^m(V_k)$.

When calculating grothendieck groups it is enough to have a filtration of comodules, hence 
the virtual Clebsch-Gordan formula may have interest for such calculations.

\section{On the  "complete reduciblilty" property  for  $\SLZ$}

Let in this section $p \in \Z $ be a prime number and let $\Z_{(p)}$ be the localization of $\Z$ at $p$. Let $m$ be an integer and let $\Z[\frac{1}{m}]$ be the localization of $\Z$ at the set $S:=\{m^i\}$.
We prove that the group schemes $\SLZ, \SL_{2,\Z_{(p)}}$ and $\SL_{2, \Z[\frac{1}{m}]}$ do not satisfy the "complete reducibility" property. Over a field of $k$ characteristic zero, any finite dimensional comodule on $\SL(2,k)$
satisfies this property.  In order to study torsion free comodules  on $\SLZ$ we introduce the notion "good filtration" and discuss the problem of classifying
torsion free finite rank comodules with a good filtration. We give explicit non trivial examples of torsion free finite rank module with a good filtration, hence many nontrivial torsion free comodules has such a filtration, and the notion may therefore have general interest in the study of comodules on a reductive group scheme over a Dedekind domain. One may define the notion "good filtration" for torsion free comodules on $\SLA$ for any Dedekind domain $A$.


Let us first introduce some notions. Let $A$ be a commutative unital ring and let $U:=\Spec(A)$. Let $G/U$ be an affine group scheme over $U$. 

\begin{definition} We say $G$ is \emph{completely reducible} iff for any surjective map of $G$-comodules $p: V \rightarrow W$ where $V$ and $W$ are finitely generated $A$-modules, there is a map
of $G$-comodules $s: W \rightarrow V$ with $ s \circ p = Id_W$.
\end{definition}

Note: A surjective map of $G$-comodules $p: V \rightarrow W$ has a $G$-comodule section $s: W \rightarrow V$ iff there is an isomorphism of $G$-comodules $V \cong K \oplus W$ where $K:=ker(p)$.
We also say $G$ satisfies the \emph{complete reducibility property} if $G$ is completely reducible.

\begin{theorem}  \label{ss} Let $k$ be a field of characteristic zero and let $G$ be an affine group scheme of finite type over $k$, with Lie algebra $\Lie(G)$. It follows $\Lie(G)$ is semi simple iff the centre of $G$ is finite and
any finite dimensional $G$-comodule $V$ is isomorphic to a direct sum $V \cong V_1 \oplus \cdots \oplus V_i$ where $V_i$ is an irreducible $G$-comodule.
\end{theorem}
\begin{proof} This is proved in \cite{DG}, Corollary 2.2, page 260.
\end{proof}

\begin{corollary} Let $k$ be a field of characteristic zero and let $G$ be an affine group scheme of finite type over $k$. If $\Lie(G)$ is semi simple it follows $G$ is completely reducible.
\end{corollary}
\begin{proof} Let $p: V \rightarrow W$ be a surjective map of $G$-comodules. It follows from Theorem \ref{ss} there are direct sum decompositions 

\[ V \cong V_1 \oplus \cdots \oplus V_i, W \cong W_1 \oplus \cdots \oplus W_j \]

where each $V_s, W_t$ is irreducible. Since for any $s,t$ it follows $\Hom_{G-comod}(V_s, W_t)$ is either zero or $k$ we must have $V \cong K \oplus W$with $K:=ker(p)$. The Corollary follows.
\end{proof}

Let $V:=A\{e_1,e_2\}$ and consider the sequence of  comodules

\[ 0 \rightarrow \Sym^{n-1}(V) \rightarrow^{f_{0,n}}  V \otimes \Sym^n(V) \rightarrow^{\pi_n} \Sym^{n+1}(V) \rightarrow 0.\]

where $g$ is the canonical multiplication map and where $f(v):=zv$ where $z:=e_1 \otimes e_2 - e_2 \otimes e_1$. It follows the sequence is an exact sequence of $\SLA$-comodules where $A$ is an arbitrary commutative unital ring.

 It follows $V \otimes \Sym^n(V)$ has a basis 

\[ V \otimes \Sym^n(V) \cong A\{ e_i \otimes e_1^n, e_i \otimes e_1^{n-1}e_2,  \cdots , e_i \otimes e_2^n\} .\]

Let for $1 \leq i \leq n-1$  

\[ v_i:= (e_1 \otimes e_2 -e_2\otimes e_1)e_1^{n-i}e_2^{i-1} :=ze_1^{n-i}e_2^{i-1}.\]

It follows we may chose a basis for $V \otimes \Sym^n(V)$ as follows:

\[ V \otimes \Sym^n(V) \cong  \]

\[  A\{ e_1 \otimes e_1^n, e_1 \otimes e_1^{n-1}e_2, ze_1^{n-1}, e_1\otimes e_1^{n-2}e_2^2, ze_1^{n-2}e_2,\ldots , e_1\otimes e_2^{n-1}, ze_2^{n-1}, e_2 \otimes e_2^n\} .\]

A section $s$ of the projection map $\pi_n: V \otimes \Sym^n(V) \rightarrow \Sym^{n+1}(V)$ would have to look as follows:

\[ s(e_2^{n+1})=e_2 \otimes e_2^n, s(e_1e_2^n)=e_1 \otimes e_2^n+ aze_2^{n-1} \]

with $a \in A$.  The element $x \in \Lie(G)$ acts as follows: $x(e_2^{n+1})=(n+1)e_1e_2^n$. Moreover $x(e_2 \otimes e_2^n)=e_1 \otimes e_2^n + ne_2 \otimes e_1e_2^{n-1}$. It follows

\[ xs(e_2^{n+1})= e_1 \otimes e_2^n + n(e_2 \otimes e_1e_2^{n-1}).\]

It follows 

\[ n(e_2 \otimes e_1e_2^{n-1}) = n e_1 \otimes e_2^n - nze_2^{n-1} \]

hence for $s(x(e_s^{n+1})=xs(e_2^{n+1})$ to hold we must have an equality

\[  nze_2^{n-1}=a(n+1)ze_2^{n-1} \]

and it follows $a:= \frac{n}{n+1} \in A$. From this it follows $\frac{1}{n+1} \in A$. Hence the projection map $\pi$ does not have a section $s$ if $A=\Z$. Hence the group scheme $\SLZ$ does not satisfy the "complete reducibility" property.
If $(p) \subseteq \Z$ is a prime ideal and if $A:=\Z_{(p)}$ is the localization of $\Z$ at $(p)$, it follows the sequence

\[ 0 \rightarrow \Sym^{p-2}(V) \rightarrow^f V \otimes \Sym^{p-1}(V) \rightarrow^g \Sym^{p}(V) \rightarrow 0\]

does not split since $p$ is not a unit in $A$. Hence $\SL_{2, \Z_{(p)}}$ does not satisfy the "complete reducibility" property.  If  $G:= \SL_{2, \Z[\frac{1}{m}]}$ there is always an integer $n$ where the projection map $\pi_n$ does not split,
hence $G$ does not satisfy the "complete reducibility" property.

\begin{theorem} \label{completereducibility} Let $m$ be any positive integer and let $p\in \Z$ be a prime number. The group schemes $\SLZ, \SL_{2, \Z_{(p)}}$ and $ \SL_{2, \Z[\frac{1}{m}]}$ do not satifsy the complete reducibility property.
\end{theorem}
\begin{proof} The proof is given above.
\end{proof}

Let $\Delta_V: V \rightarrow R \otimes_{\Z} V$ be a comodule that is finitely generated as $\Z$-module and let $T(V) \subseteq V$ be the torsion sub module. There is a direct sum decomposition $V \cong E(V) \oplus T(V)$
where $E(V) \cong \Z^n$ is a free $\Z$-module of rank $n$. Let $p \in \Z$ be a prime and let $T(V)_p$ be the set of elements $v$  in $T(V)$ with $p^mv=0$ for some integer $m\geq 1$. It follows there is a direct sum decomposition

\[ T(V) \cong \oplus_{p \in J} T(V)_p \]

where $J:=\{ p_1,..,p_l\} \subseteq \Z$ is a finite set of distinct prime numbers.

\begin{proposition} \label{notirreducible} Let $I:=(m) \subseteq \Z$ be an ideal. Let $W:=\Sym^*(V)/I\Sym^*(V)$ and let $W_n:= \Sym^n(V)/I\Sym^n(V)$. There is a canonical comodule structure $ \Delta_I: W \rightarrow R \otimes W$ with the property
that the canonical projection map $\pi: \Sym^*(V) \rightarrow W$ is a map of comodules. For each integer $n \geq 1$ we get a comodule $\Delta_{I,n}: W_n \rightarrow R \otimes W_n$ and there is an exact sequence of comodules

\[ 0 \rightarrow \Sym^n(V) \rightarrow^\phi \Sym^n(V) \rightarrow^p W_n \rightarrow 0 \]

where $\phi(v):=mv$. It follows $Im(\phi) \subsetneq \Sym^n(V)$ is a strict sub-comodule and the comodule $\Sym^n(V)$ is not irreducible. Assume $V$ is a comodule that is finitely generated as $\Z$-module with a direct sum decomposition
$V \cong \Z^n \oplus T(V)$ where $T(V)$ is the torsion submodule. There is an exact sequence of comodules

\[  0 \rightarrow T(V) \rightarrow V \rightarrow V/T(V) \rightarrow 0. \]

For any prime $p \in J$ it follows $T(V)_p \subseteq T(V)$ is a sub comodule. It follows the direct sum $T(V) \cong \oplus_{p \in J}T(V)_p$ is a direct sum of comodules. Let $E(V):=V/T(V)$. It follows
$V \cong E(V)\oplus T(V)$ is a direct sum of comodules.

\end{proposition}
\begin{proof} The proof is an exercise.
\end{proof}

If $k$ is a field of characteristic zero and if $V$ is any finite dimensional $\SL(2,k)$-module, it follows there is a direct sum decomposition $V \cong \oplus_{i \in I} \Sym^{d_i}(V)$ of representations. The symmetric powers $\Sym^d(V)$ are irreducible
when $V:=k\{e_1,e_2\}$ is the standard $\SL(2,k)$-module.  Proposition \ref{notirreducible} proves this property does not hold for the group scheme $\SLZ$. The symmetric power of the standard module $\Sym^d(V)$ is not irreducible
and there is no decomposition of the tensor product $\Sym^n(V) \otimes \Sym^m(V)$ into a direct sum of irreducible modules. By the above example it follows $V \otimes \Sym^n(V)$ is not isomorphic to the direct sum $\Sym^{n-1}(V) \oplus \Sym^{n+1}(V)$
as comodules. Hence the "classification" of $\SLZ$-modules of finite rank over $\Z$ is more complicated than the case of $\SL(2,k)$.

Proposition \ref{notirreducible} proves that for any comodule $V$ that is finitely generated as $\Z$-module, if we want to "classify" $V$, we may "classify" the torsion free part $E(V)$ and the $p$-part $T(V)_p$ separately, since there is a direct sum of comodules

\[  V \cong E(V) \oplus (\oplus_{p \in J} T(V)_p) .\]

\begin{definition} Let $V$ be the standard comodule on $\SLZ$ from Definition \ref{standard} and let $W$ be a torsion free $\SLZ$-comodule which is of finite rank as $\Z$-module. We say $W$ has a \emph{good filtration} if there is a finite filtration of comodules  $F_i \subseteq W$  with 

\[  F_i/F_{i+1} \cong \Sym^{d_i}(V) \]

with $d_i \geq 0$ an integer for all $i$. 

\end{definition}

By Corollary \ref{CGfiltration} the following holds: For any pair of integers $1 \leq n \leq m$ it follows the tensor product $\Sym^n(V)\otimes \Sym^m(V)$ has a good filtration. We may ask for a classification
of torsion free comodules on $\SLZ$ having a good filtration.

\begin{example}  \label{notirreducible} An elementary example. \end{example}

Let $G:=\SLZ$ and let $V:=\Z\{e_1,e_2\}$ be the standard $G$-comodule from Definition \ref{standard}. Define for any integer $m \in \Z$ a map $f_m: V \rightarrow V$ by $f_m(v):=mv$.
Let $I:=(m) \subseteq \Z$ be the ideal generated by $m$ and assume $I \neq (0),(1)$. Consider the following exact sequence pf $\Z$-modules:

\begin{align}
&\label{exct}  0 \rightarrow V \rightarrow^{f_m} V \rightarrow  Q_m \rightarrow 0,
\end{align}

with $Q_m \cong \Z/(m)\{e_1,e_2\} \neq (0)$. It follows \ref{exct} is an exact sequence of $G$-comodules, hence $Im(f_m) \subsetneq V$ is a strict sub $G$-comodule. It follows $V$ has an infinite set of strict sub comodules hence $V$ is not an irreducible comodule.

Note: Let $k$ be the field of rational numbers and let $G_k:=\SL(2,k)$. By \cite{serre} it follows the representation rings $\operatorname{R}(G) \cong \operatorname{R}(G_k)$ are isomorphic. One wants to interpret
the results in this section in terms of this isomorphism. 

Let $R:=A[x_{ij}]/(D-1)$ with $\SLA:=\Spec(R)$ and let $T: R \rightarrow R$ be the following map:

\[ T(x_{11}):=x_{11}, T(x_{12}):=x_{21}, T(x_{21}):=x_{12}, T(x_{22}):=x_{22}.\]

Hence $T^2=Id$ is the idenitity and the map $T$ corresponds to taking the tranpose of a matrix. It follows $T$ is an automorphisms of $A$-algebras. 
Recall the inversion map $S: R \rightarrow R$:

\[ S(x_{11}):= x_{22}, S(x_{12}):=-x_{12}, S(x_{21}):=-x_{21}, S(x_{22}):=x_{11}.\]

It follows $S$ is an automorphism of $A$-algebras. Let $\T:= T \circ S$. The following holds:

\begin{lemma} There is an equality of maps $\T = T \circ S = S \circ T$. There is an equality $\T \otimes \T \circ \Delta_R=  \Delta_R \circ \T$. Moreover $\T \circ \T = Id$.
\end{lemma}
\begin{proof} As an example: There is an equality

\[ \T \otimes \T (\Delta_R(x_{11}))= x_{22}\otimes x_{22} + x_{21} \otimes x_{12} = \Delta_R(x_{22})= \Delta_R(\T(x_{11}) .\]

The Lemma follows.
\end{proof}

We are going to study the classical dual of a classical representation of an abstract group.

\begin{example} The "classical" dual of a "classical" representation. \end{example}

If $G$ is an abstract group acting from the left on a finite dimensional vector space $V$ over a field $k$, with $\sigma: G \times V \rightarrow V$ the action, we may dualize the action 
$\sigma$ as follows. Define the \emph{classical dual action} 

\[\sigma^*: G \times V^* \rightarrow V^* \]

by $\sigma^*(g, \psi)(v):= \psi(g^{-1}v)$. There is a canonical action of $G$ on the direct sum $V \oplus V^*$ defined by $g(v, \psi):= (gv, g\psi)$ and the trace map $tr: V \oplus V^* \rightarrow k$
is invariant with respect to this action. Moreover the canonical map $\phi: V \rightarrow (V^*)^*$ is an isomorphism of $G$-modules. We want to do something similar for comodules.

\begin{example} The "tranpose"  of a free comodule. \end{example}

We can try to "mimic" the construction of the dual action of a representation of an abstract group and try to do something similar for a representation of a group scheme. 
For a free comodule $(W , \Delta)$ if finite rank over $A$, there will in general  be several comodule structures on the dual $W^*$. Now we will define the \emph{transpose comodule} $\D^{\T}$ of $W$.


Let $W:=A\{e_1,..,e_n\}$ be the free $A$-module of rank $n$ on the elements $e_i$, and let $\Delta: W \rightarrow R \otimes W$ be a comodule structure with $\Delta(e_i):= \sum_j a(i)_j \otimes e_j$ for $i=1,..,n$.

 Let $W^*:=A\{x_1,..,x_n\}$ with $x_i:=e_i^*$ and define the map 

\[ \Delta_1: W^* \rightarrow R \otimes W^* \]

by 

\[  \Delta_1(x_i):= \sum_j a(i)_j \otimes x_i .\]

Define $\Delta^{\T}: W^* \rightarrow R \otimes W^*$ by $\Delta^{\T}:= \T \otimes 1 \circ \Delta_1$.

\begin{proposition}\label{dual}  The maps $\Delta_1, \Delta^{\T}$ are comodule structures on $W^*$. The canonial isomorphism of $A$-modules
$\phi: W \rightarrow (W^*)^*$ is an isomorphism of comodules. 
\end{proposition} 
\begin{proof} One checks $\Delta_1$ is a comodule structure on $W^*$. Assume

\[ \Delta_1(x_i)= \sum_l a(i)_l \otimes x_l .\]

It follows 

\[ 1 \otimes \Delta_1(\Delta_1(x_i))= \sum_{l,k} a(i)_l \otimes a(l)_k x_k .\]

Assume 

\[ \Delta_R(a(i)_l)= \sum_p u(i,l)_p \otimes v(i,l)_p.\]

It follows there is an equality

\[ 1 \otimes \Delta_1(\Delta_1(x_i))= \sum_{l,k} a(i)_l \otimes a(l)_k \otimes x_k =\]

\[ \sum_{l,p} u(i,l)_p \otimes v(i,l)_p \otimes x_l .\]

It follows

\[ 1\otimes \Delta^{\T}(\Delta^{\T}(x_i))= \sum_{l,k} \T(a(i)_l) \otimes \T(a(l)_k) \otimes x_k =\]

\[  \sum_{l,p} \T(u(i,l)_p) \otimes \T(v(i,l)_p) \otimes x_l =\]

\[ \T \otimes \T \otimes 1 (\Delta_R\otimes 1(\Delta_1(x_i))= \Delta_R\otimes 1(\T \otimes 1(\Delta_1(x_i))= \Delta_R\otimes 1(\Delta_2(x_i)) \]

hence $\Delta_R \otimes 1 \circ \Delta^{\T} = 1 \otimes \Delta^{\T} \circ \Delta^{\T}$. The rest of the proof is similar and the Proposition is proved.
\end{proof}

\begin{definition} Let $(W, \Delta)$ be a comodule with $W$ a free $A$-module of rak $n$.
 We define $(W^*, \Delta^{\T} )$ from Proposition \ref{dual} to be the \emph{transpose comodule} of $(W, \Delta)$. 
\end{definition}

\begin{example} The transpose of a comodule on $\SLZ$. \end{example}

As an example: If $G:=\SLZ$  it follow any comodule $W$ that is finitely generated as $\Z$-module, is the direct sum (a direct sum of comodules) $W \cong E(W) \oplus T(W)$ with $E(W) \cong \Z^n$ a free rank $n$ module on $\Z$. The module $E(W)$ is the torsion free part of $W$. We may for such modules always define the transpose comodule $W^* \cong E(W)^*$ since $T(W)^* =(0)$ is zero. Hence for any comodule $W$ on $\SLZ$ that is finitely generated as $\Z$-module we may always define the transpose comodule $W^*$ using Proposition \ref{dual}. 




\begin{example} A class of comodules with a good filtration. \end{example}

Let $k$ be the field of rational numbers and let $\SL(2,k)$ be the pull back of $\SLZ$ to $k$. There is by \cite{serre} an isomorphism

\[ \eta: \operatorname{K}(\SLZ) \cong  \operatorname{K}(\SL(2,k)) \]

of abelian groups and $ \operatorname{K}(\SL(2,k))$ is generated by the classes $[\Sym^d(V\otimes k)]$ where $V\otimes k$ is the standard comodule on $\SL(2,k)$. The comodule $V\otimes k$ is the pull back of $V$ to $\SL(2,k)$.

\begin{lemma} \label{lemmagr}The group $ \operatorname{K}(\SLZ)$ is generated by the classes $[\Sym^d(V)]$ where $d \geq 0$ is an integer and $(V, \Delta)$ is the standard comodule. There is an equality of classes

\[  [\Sym^d(V)^* ] = [\Sym^d(V)] \]
for any integer $d \geq 0$.
\end{lemma}
\begin{proof} Let $ x \in \operatorname{K}(\SLZ)$ be any class. It follows  $\eta(x)=\sum_i n_i[\Sym^{d_i}(V\otimes k)]$ for some integers $d_i \geq 0$ and $n_i \in \Z$. It follows $y:=\sum_i n_i[\Sym^{d_i}(V)] \in \operatorname{K}(\SLZ)$
and since $\eta$ is injective and $\eta(x-y)=0$, it follows $x=y$. The second statement follows using a simlar argument since there is an isomorphism 

\[ \Sym^d(V)^* \otimes_{\Z} k \cong \Sym^d(V \otimes_{\Z} k) \]

of representations of $\SL(2,k)$.
\end{proof}

Note: Lemma \ref{lemmagr} and Corollary \ref{CGfiltration} suggests that many torsion free comodules on $\SLZ$ have a good filtration. Since the group scheme $\SLA$ fails to be completely reducible for a general Dedekind domain,
the existence of a good filtration for a large class of torsion free comodules may have interesting consequences for this study. In the following we give more nontrivial examples of torsion free comodules with a good filtration.

Let $\phi: V \rightarrow V$ be defined by $\phi(e_1):=e_2, \phi(e_2):=-e_1$. Let $\Sym^n(V):=\Z\{v_0,..,v_n\}$ with $v_i:=e_1^{n-i}e_2^i$ and let $x_i:=v_i^*$ be the dual basis. Let $\Sym^n(\phi):=\phi^n$ and $\Sym^*(\phi):=\phi^*$.

Let 

\[  \T\otimes \phi^*: R \otimes \Sym^*(V) \rightarrow R \otimes \Sym^*(V) \]

be defined by $\T \otimes \phi^*(x \otimes v):=\T(x)\otimes \phi^*(v)$. Let $(V, \Delta)$ be the standard comodule, let $(\Sym^n(V), \Delta_n)$ be the $n$'th symmetric product and let $\Sym^n(V)^*, \Delta^{\T}_n)$ be the tranpose as defined in Proposition \ref{dual}.

It follows 

\[ \T \otimes \phi^*(\Delta(e_1))= \Delta(e_2), \T\otimes \phi^*(\Delta(e_2))=-\Delta(e_1). \]

It follows 

\[      \T \otimes \phi^*(\Delta_n(v_i))=           \T \otimes \phi^*(\Delta_n(e_1^{n-i}e_2^i))=(-1)^i \Delta_n(e_1^ie_2^{n-i})=(-1)^i\Delta_n(v_{n-i}) .\]

Define the elements $a(i)_l \in R$ as follows: 

\[ \Delta_n(v_i):=a(i)_0 \otimes v_0+ a(i)_1 \otimes v_1 + \cdots + a(i)_n \otimes v_n.\]

The elements are well defined since $v_0,..,v_n$ is a basis for $R\otimes \Sym^n(V)$ as $R$-module. The elements $a(i)_l$ are "symmetric" with respect to the involution $\T$ in a sense expressed in the following Lemma:

\begin{lemma}\label{symmetry}  For all $i,l=0,1,..,n$ there is an equality

\[    \T(a(i)_l)=(-1)^{i+l}a(n-i)_{n-l}.\]

\end{lemma}
\begin{proof}
It follows 

\[  \T \otimes \phi^*(\Delta_n(v_i))= \]

\[ \T(a(i)_0) \otimes v_n + \cdots + (-1)^l \T(a(i)_l) \otimes v_{n-l} + \cdots + (-1)^n \T(a(i)_n)\otimes v_0 =\]

\[ (-1)^n \T(a(i)_n)\otimes v_0 + \cdots + (-1)^l\T(a(i)_l) \otimes v_{n-l} + \cdots + \T(a(i)_0)\otimes v_n =\]

\[ (-1)^i \Delta_n(v_{n-i})=   \]

\[ (-1)^i a(n-i)_0\otimes v_0 + \cdots + (-1)^i a(n-i)_{n-l}\otimes v_{n-l} + \cdots + (-1)^i a(n-i)_n \otimes v_n.\]

It follows there is for all $i,l=0,..,n$ an equality

\[    \T(a(i)_l)=(-1)^{i+l}a(n-i)_{n-l}.\]

The Lemma follows.
\end{proof}

\begin{theorem} \label{main} Let $(V ,\Delta)$ be the standard $\SLZ$-comodule and let $W_1,W_2$ be $\SLZ$-comodules on the form $\Sym^n(V), \Sym^m(V^*)$ or $(\Sym^k(V)^*, \D^{\T}_k)$ for integers $n,m,k \geq 0$. It follows the tensor power $W_1 \otimes W_2$ has a good filtration.
\end{theorem}
\begin{proof} There is an isomorphism of comodules $V \cong V^*$ hence $\Sym^n(V) \cong \Sym^n(V^*)$ are isomorphic as comodules. If $v_i:=e_1^{n-i}e_2^i, x_i:=v_i^*$ are bases for $\Sym^n(V)$ and $\Sym^n(V)^*$, there is an isomorphism
of $\Z$-modules

\[ f: \Sym^n(V) \cong \Sym^n(V)^* \]

defined by $f(v_i):=(-1)^ix_{n-i}$. By definition

\[ \Delta_n(v_i)= \sum_l a(i)_l \otimes v_l \]

and

\[ \Delta^{\T}_n(x_i):= \sum_l \T(a(i)_l) \otimes x_l .\]

We get

\[ 1\otimes f(\Delta_n(v_{n-i}))= a(n-i)_0 \otimes x_n - a(n-i)_1 \otimes x_{n-1} + \cdots + \]

\[   (-1)^{n-l} a(n-i)_{n-l} \otimes x_l + \cdots + (-1)^n a(n-i)_n \otimes x_0.\]

We also get

\[ \Delta^{\T}_n(f(v_{n-i}))= (-1)^{n-i}\Delta^{\T}_n(x_i)= \]

\[ (-1)^{n-i} \T(a(i)_0)\otimes x_0 +  (-1)^{n-i}\T(a(i)_1)\otimes x_1+  \cdots +  (-1)^{n-i}   \T(a(i)_l)\otimes x_l + \cdots + \]

\[ (-1)^{n-i}   \T(a(i)_n)\otimes x_n .\]





By Lemma \ref{symmetry} there is an equality  $\T(a(i)_l)=(-1)^{i+l}a(n-i)_{n-l}$. Hence $\Delta^{\T}_n(f(v_i))= 1\otimes f(\Delta_n(v_i))$ and we have proved $f$ is an isomorphism of comodules. The Theorem follows from Lemma \ref{CGfiltration}.

\end{proof}

Note: If $k$ is the field of rational numbers and if $V(k):=V\otimes k$ is the standard comodule on $\SL(2,k)$, there is an isomorphism $V(k) \cong V(k)^*$ and $\Sym^n(V(k)) \cong \Sym^n(V(k)^*)$. There is a "canonical" isomorphism
of comodules $\tau: \Sym^n(V(k)^*) \cong \Sym^n(V(k))^*$ defined by  

\[ \tau(\phi_1 \cdots \phi_n)(v_1\cdots v_n):= \sum_{\sigma} \phi_{\sigma^{-1}(1)}(v_1) \cdots \phi_{\sigma^{-1}(n)}(v_n) ,\]

hence we get a "canonical"  isomorphism $\tau: \Sym^n(V(k)) \cong  \Sym^n(V(k))^*$.

The isomorphism  in Theorem \ref{main} differs from $\tau$. There is a map similar to $\tau$ defined for $\SLZ$ but it is no longer an isomorphism of comodules. The isomorphism $\tau$ is constructed in \cite{fulton},
in section B3, page 476.

\section{Generalities on duals and left and right comodules on $\SLZ$}

I this section we include some general results on left and right comodules on  Hopf algebras. We also give an explicit construction of the dual $(V^*, \Delta^*)$ of any comodule $(V, \Delta)$ on a Hopf algebra $R$ in the case when $V$
is a free and finite rank module over the base ring $k$. For the group scheme $\SLZ$ any torsion free finitely generated comodule is free, hence we get a functorial definition of the dual of any torsion free and finitely generated
 comodule on $\SLZ$. Some of the proofs in the section are left to the reader as an exercise. They are standard proofs involving commutative diagrams and a bit of homological algebra.

We use the construction of the dual to give an explicit and elementary construction of a torsion free comodule $W$ on $\SLZ^{op}$ where $W$ and the dual $W^*$ have the same weight space decomposition for the induced action of $\Lie(\SLZ^{op})$, but where $W$ and $W^*$ are non isomorphic as comodules on $\SLZ^{op}$. Hence a torsion free comodule is not determined by its weight space decomposition as is the case over a field $k$ of characteristic zero.

Note the following: Let $R$ he a Hopf algebra in the sense of Jantzens book. An affine group scheme is a functor

\[  h_R: \alg \rightarrow \grp \]

where $\alg$ is the category of commutative $A$-algebras and $\grp$ is the category of groups. It is defined by $h_R(B):=\Hom_{A-alg}(R,B)$. Let $V \cong A\{e_1,..,e_n\}$ be a free rank $n$ module on $A$ and define the following functor:

\[ \GL(V): \alg \rightarrow \grp \]

by $\GL(V)(B):=\GL_B(B \otimes_A V)$.  A "representation" in the language of Jantzen is a natural transformation of functors $\eta: h_R \rightarrow \GL(V)$. Given such a natural transformation $\eta$ we may dualize $\eta$ in the following way:
For any $A$-algebra $B$ and any element $g \in h_R(B)$ we get an element $\eta_B(g)^{-1} \in \GL_B(B\otimes_A V) $ and we may for any element  $f \in (B\otimes_A V)^* \cong V^*\otimes_A B$ define $\eta^*_B(g)(f):= f \circ \eta_B(g)^{-1}$.
We get the "dual representation"

\[ \eta^*: h_R \rightarrow \GL(V^*) ,\]

which is proved to be a natural transformation of functors. When doing this we do not "worry" about left and right representations and if one wants to define the dual directly on the comodule $(V,\D)$ corresponding to $V$ one has to do this.
The issue of left and right comodules is not mentioned in the books of Demazure, Gabriel and Jantzens and here I will discuss this construction in detail.

 Note: The category $\alg$ of all commutative $A$-algebras is not a small category, hence if you have issues with non standard set theory, you may prefer the approach using comodules. Else you can use the definition given in the books of of 
Demazure, Gabriel and Jantzen.
The approach in this section gives an explicit construction of the dual comodule $(V^*, \D^*)$ of any comodule $(V,\D)$ on $\SLZ$ where $V$ is a finitely generated and torsion free $\Z$-module. If you take the "derivative" of the dual action $\eta^*$,
you get the dual action of the Lie algebra $\mathfrak{sl}(n,\Z)$ on $V^*$, but there is a strict inclusion $\operatorname{U}(\mathfrak{sl}(2,\Z)) \subsetneq \Dist(\SLZ)$. You may not use the Lie algebra or the enveloping algebra
to study the isomorphism class of the dual comodule $(V^*, \D^*)$. Two comodules $V,W$ on $\SLZ$ may be isomorphic as representations of the Lie algebra, but non-isomorphic as comodules.

Let in this section $\SLZ:=\Spec(R)$ where $(R, \D_R, S, \e)$ is the Hopf algebra introduced in section 1. Let $ \sigma: R\otimes_{\Z} R \rightarrow R \otimes_{\Z} R$ be defined by $\sigma(a\otimes b):= b \otimes a$ and let $\tD_R:= \sigma \circ \D_R$.
It follows 

\[ \Delta_R(x_{11}):= x_{11} \otimes x_{11} + x_{12} \otimes x_{21}, \Delta_R(x_{12}):=x_{11} \otimes x_{12} + x_{12} \otimes x_{22} ,\]

\[ \Delta_R(x_{21}):= x_{21} \otimes x_{11} + x_{22} \otimes x_{21}, \Delta_R(x_{22}):= x_{21} \otimes x_{12} + x_{22} \otimes x_{22}.\]

and

\[ \tD_R(x_{11}):= x_{11} \otimes x_{11} + x_{21} \otimes x_{12}, \tD_R(x_{12}):=x_{12} \otimes x_{11} + x_{22} \otimes x_{12} ,\]

\[ \tD_R(x_{21}):= x_{11} \otimes x_{21} + x_{21} \otimes x_{22}, \tD_R(x_{22}):= x_{12} \otimes x_{21} + x_{22} \otimes x_{22}.\]

One checks $(R, \tD_R, S,\e)$ is a Hopf algebra.

\begin{definition} The Hopf algebra $(R, \tD_R, S, \e)$ is the \emph{opposite Hopf algebra} of $(R, \D_R, S, \e)$.
\end{definition}

\begin{example} \label{dualgroup} A representation and its dual. \end{example}

Let $V:=\Z\{e_1,e_2\}$ and let $B$ be a commutative $\Z$-algebra. Let $V \otimes_{\Z}B \cong \{e_1 \otimes 1, e_2\otimes 1\}B$ and let $C:=\{e_1\otimes 1, e_2\otimes 1\}$.  Let $V^*:=\Z\{x_1,x_2\}$ and let $B\otimes_{\Z} V^*:=B\{1 \otimes x_1, 1 \otimes x_2\}$. Let $C^*:=\{1 \otimes x_1, 1\otimes x_2\}$.

If $u:=e_1 \otimes u_1, + e_2\otimes u_2$ with $u_i \in B$ write

\[
[u ]_C:=  
\begin{pmatrix}   u_1  \\
   u_2
\end{pmatrix}.
\]

Hence we let $[u]_C$ be a "column vector". If $f:=v_1x_1 + v_2 x_2$ we write $[f]_{C^*}:=[v_1,v_2]$. Hence we let $[f]_{C^*}$ be a row vector. There is a pairing $tr: V^* \oplus V \rightarrow \Z$ and we have

\[ tr(f, u):=f(u):= v_1u_2+v_2u_2.\]

There are two functors 

\[  h_R, \GL(V): \zalg \rightarrow \grp \]

defined by $h_R(B):=\Hom_{\Z-alg}(R, B)$ and $\GL(V)(B):=\GL_B(V\otimes_{\Z} B)$. There is a natural tranformation $\rho: h_R \rightarrow \GL(V)$ defined as follows:

\[  \rho_B: h_R(B) \rightarrow \GL_B(V \otimes_{\Z} B) \]

By definition $V \otimes B \cong B^2$. An element $g \in h_R(B)$ is by definition a $2 \times 2$ matrix $g:=(a_{ij})$ with $det(g)=1$. Hence we may let $g$ act on $V \otimes B$ using matrix multiplication from the left. One checks $\rho$ defines a natural
transformation of functors. The natural transformation $\rho$ is by definition a representation of $\SLZ$ on the module $V$. We may calculate its corresponding comodule and one finds that the  comodule $\D_V: V \rightarrow V \otimes R$ is the following:

\[ \Delta_V(e_1):= e_1  \otimes x_{11}+ e_2 \otimes x_{21}, \D_V(e_2):= e_1 \otimes x_{12} + e_2 \otimes x_{22}.\]

We may dualize $\rho$ to a natural transformation

\[ \rho^*:  h_R \rightarrow \GL(V^*) \]

and by definition  $\rho_B^*: h_R(B) \rightarrow \GL_B(B \otimes V^*)$ is defined as follows: Given any element $f \in B\otimes V^* \cong (V \otimes B)^*$ we let $\rho^*_B(g)(f):= f \circ \rho_B(g^{-1})$. We get a comodule structure
$ \D_{V^*}: V^* \rightarrow R^{op} \otimes V^*$ defined as follows:

\[ \D_{V^*}(x_1):= x_{22} \otimes x_1 - x_{12} \otimes x_2, \D_{V^*}(x_2):= -x_{21} \otimes x_1 + x_{11} \otimes x_2.\]

Let $g^{-1}:=(b_{ij})$. The action $\rho_B^*(g^{-1})(f)$ is given as follows: Multiply the vector $[f]_{C^*}:=[v_1,v_2]$ from the right with the matrix $(b_{ij})$. Hence the dual $V^*$ of a left $R$-comodule $V$ is a right $R^{op}$-comodule.

\begin{example}\label{involution}  A involution of an abstract group \end{example}

If $G$ is an abstract group with multiplication $m: G \times G \rightarrow G$ defined by $m(g,h):=gh$, define the following operation $m^*: G \times G \rightarrow G$ by $m^*(g,h):=g \bullet h:=hg$. It follows for any $g,h,k \in G$ that

\[  (g \bullet h ) \bullet k = k(hg)=(kh)g = g \bullet (h \bullet k) \]
hence the operation $\bullet $ is associative. It follows $(G, \bullet)$ is a group and the map $\phi: (G, m) \rightarrow (G, \bullet)$ defined by $\phi(g):= g^{-1}$ is an isomorphim of groups with the property that $\phi^2 = Id_G$. 


We want to define something similar for group schemes.

\begin{lemma} The 4-tuple $(R, \tD_R, S, \e)$ is a Hopf algebra. The map $S: R \rightarrow R$  sits in a commutative diagram of map of commutative rings

\[
\diagram       R \dto^{\tD_R}    \rto^S  &   R \dto^{\D_R}  \\
                        R\otimes R \rto^{S\otimes S} & R \otimes R.
\enddiagram
\]
There is also an equality of maps $\tD_R \circ S = S\otimes S \circ \D_R$.
\end{lemma}
\begin{proof} The proof is an exercise.
\end{proof}

We may view the map $S$ as a group scheme version of the involution $\phi$ from Example \ref{involution}. The map $S$ induce a map of group schemes $\phi_S: \SLZ \rightarrow \SLZ^{op}$  where
$\SLZ^{op}:=\Spec(R)$ with $\tD_R$ as comultiplication. The map $S$ is an isomorphism of Hopf algebras hence $\phi_S$ is an isomorphism of group schemes. We may do a similar construction for any affine group scheme $H$.

\begin{definition}
Let  $H$ be a group scheme. Let $\co^H$ denote the category of comodules $\Delta: V \rightarrow R \otimes V$ with the obvious maps.  Let $\co_H$ denote the category of maps of comodules $\Delta: V \rightarrow V \otimes R$
with the obvious maps. We let $\co^H$ be the \emph{category of right $H$-comodules} and we let $\co_H$ be the \emph{category of left $H$-comodules}.
\end{definition}

Define for any $\Z$-module $V$ the following maps: $\tau_1: R \otimes V \rightarrow V \otimes R$ by $\tau_1(x \otimes v):= v \otimes x$. Define $\tau_2: V \otimes R \rightarrow R \otimes V$ by $\tau_2(v\otimes x):= x \otimes v$. It follows $\tau_1$ and $\tau_2$ are inverses of each other.

Let $\D_V: V \rightarrow R \otimes V$ be an $H$-comodule and let $\D_W: W \rightarrow W \otimes R$ be an $H^{op}$-comodule. Let $\D_1:= \tau_1 \circ \D_V$ and let $\D_2:= \tau_2 \circ \D_W$.

Define the maps $F,G$ as follows:

\[ F(V, \D_V):= (V, \D_1),  \]

where 

\[ \D_1: V \rightarrow V \otimes R \]

and

\[  G(W, \D_W):=(W, \D_2) \]

where

\[ \D_2: W \rightarrow R \otimes W.\]

There is a commutative diagram

\[
\diagram       V  \rto^{\D_V} \dto^{\D_V}   &    R \otimes V \dto^{1\otimes \D_V}   \\
                       R \otimes V\rto^{\D_R \otimes 1}  &  R \otimes R \otimes V.
\enddiagram
\]

Assume  for any  $v \in V$ we may write

\[ \D_V(v):= \sum_i a_i \otimes v_i ,\]

\[ \D_V(v_i):= \sum_j  b(i)_j \otimes v(i)_j ,\]

and

\[ \D_R(a_i):= \sum_k u(i)_k \otimes w(i)_k .\]

Since the diagram commutes there is an equality

\[ \sum_{i,j} a_i \otimes u(i)_j \otimes v(i)_j = \sum_{i,k} u(i)_k \otimes w(i)_k \otimes v_i.\]

We want to prove there is a commutative diagram

\[
\diagram       V  \rto^{\D1} \dto^{\D1}   &     V \otimes R \dto^{ \D_1 \otimes 1}   \\
                       R \otimes V\rto^{1\otimes \tD_R  }  &   V\otimes R \otimes R.
\enddiagram
\]

There is a canonical isomorphism $\eta: R\otimes R \otimes V \cong V \otimes R \otimes R$ defined by $\eta(x\otimes y \otimes v):=v \otimes y \otimes x$.

We get

\[ \D_1 \otimes 1(\D_1(v))= \sum_{i,j} v(i)_j \otimes b(i)_j \otimes a_i\]

and

\[ 1\otimes \tD_R(\D_1(v))= \sum_{i,k} u(i)_k \otimes w(i)_k \otimes v_i ,\]

and since $\eta$ is an isomorphism the claim follows.

\begin{lemma} \label{involution}  The maps $F,G$ defined above give rise to functors

\[ F : \co^H \rightarrow \co_{H^{op}}\text{ and  } G: \co_{H^{op}} \rightarrow \co^H .\]
There are isomorphisms of functors $F \circ G =Id, G \circ F \cong Id$ hence $F$ and $G$ are equivalences of categories.
\end{lemma}
\begin{proof} Using methods similar to the above calculation one checks for any $H$-comodule $(V, \D_V)$ the pair $(V, \D_1)$ is a comodule on $H^{op}$. There is for any map of comodules $\phi$ an induced map $F(\phi)$ of $H^{op}$-comodules and one checks $F$ and $G$ are well defined functors.
Given any $H$-comodule $(V, \D_V)$ it follows  $(G\circ F)(V, \D_V))$ is the comodule $(V, \D_3)$ where $\D_3:= \tau_2 \circ \tau_1 \circ \D_V =\D_V$, hence $(G \circ F)(V , \D_V)=(V ,\D_V)$. The rest of the Lemma is an exercise.

\end{proof}

\begin{lemma} \label{iso} Let $k$ be any commutative unital ring and let $V$ be a free finite rank $k$-module. Let $V^*:=\Hom_k(V,k)$. There is an isomorphism of left $R$-modules

\[ \rho_1 : R\otimes_k V^* \cong \Hom_R(V\otimes_k R, R) \]

defined by $\rho_1(a\otimes f)(v\otimes b):=af(vb)$, where $a\otimes f \in R\otimes_k V^*, v\otimes b \in V\otimes_k R$. There is an isomorphism of right $R$-modules

\[ \rho_2: V^*\otimes_k R \cong \Hom_R(R\otimes_k V,R) \]

defined by $\rho_2(f\otimes a)(b \otimes v):= f(a(bv))$ where $f\otimes a\in V^*\otimes_k R, b\otimes v\in R\otimes_k V$. Let $\psi: V\otimes_k R \cong (V\otimes_k R)^{**}$ be the canonical iomorphism.
We get a canonical isomorphism 

\[  s_1:= \rho_2^{-1} \circ \rho_1^* \circ \psi: V\otimes_k R \cong V^{**}\otimes_k R, \]

 and there is an equality $s_1 \cong t \otimes 1$ with $t:V \cong V^{**}$ the canonical isomorphism. We also get a canonical isomorphim $ s_2: R\otimes_k V \cong (R\otimes_k V)^{**} \cong R\otimes_k V^{**}$ and there is an equality
$ s_2 \cong 1 \otimes t$.
\end{lemma}
\begin{proof} The proof is an exercise.
\end{proof}

Let $\D_V: V \rightarrow V \otimes_k R$ where $(R, \D_R, S, \e)$ is a Hopf algebra over $k$ and $V$ is a free finite rank $k$-module. We get an induced map

\[\phi: V\otimes_k R \rightarrow V\otimes_k R\]
and an induced dual map

\[ \phi^*: (V\otimes_k R)^* \rightarrow (V\otimes_k R)^* \]

and an induced map  $  \D_{V,1}: V^* \rightarrow (V\otimes_k R)^* \cong R \otimes_k V^*$

defined by $\D_{V,1}(f):= \phi^*(1 \otimes f)$ where $\phi^*: V^*\otimes_k R \rightarrow V^*\otimes_k R$ is the induced map from the isomorphism $\rho_2$ from Lemma \ref{iso}. 

\begin{proposition} \label{functor1} Let $H:=\Spec(R)$ and define $F(V, \D_V):=(V^*, \D_{V,1})$. It follows $F$ defines a functor

\[  F: \co_H \rightarrow \co^H .\]
\end{proposition}
\begin{proof} The pair $(V, \D_V)$ is a comodule and this implies the map $\D_V$ form two commutative diagrams. One proves this implies the map $\D_{V,1}$ satisfies the  commutative diagrams needed
for it to be a comodule in $\co^H$. Hence it follows $(V^*, \D_{V,1})$ is an object in $\co^H$. One moreover proves that for any map $f: (V, \D_V) \rightarrow (W, \D_W)$ in $\co_H$ it follows we get a map

\[  F(f): (W^*, \D_{W, 1}) \rightarrow (V^*, \D_{V, 1}) \]

in $\co^H$, and for two composable maps $f, g$ in $\co_H$ it follows $F(f \circ g)=F(g) \circ F(g)$. We get a well defined functor and the Proposition follows.

\end{proof}

Let $\D_V: V \rightarrow V \otimes R$ be a comodule on $R$ and consider the map

\[  \D_{V,1}: V^* \rightarrow  R\otimes V^* \]

from Proposition \ref{functor1}. Let $S$ be the involution of $R$ and consider the canonical map $S\otimes 1: R\otimes V^* \rightarrow R^{op} \otimes V^*$ where $R^{op}$ is the "opposite" Hopf algebra of $R$. We define a composed map

\[ \D_{V^*}:=S\otimes 1 \circ \D_{V,1}: V^* \rightarrow R^{op} \otimes V^*.\]

\begin{proposition} \label{functor2} Let $(V, \D_V)$ be a comodule on $R$ where $V$ is a free and finite rank $k$-module and let $H:=\Spec(R)$. Define $D^{LR}(V, \D_V):=(V^*, \D_{V^*})$. The map $D^{LR}(-)$ gives rise to  a functor

\[  D^{LR}_H: \co_H \rightarrow \co^{H^{op}} .\]

There is similarly a functor

\[ D^{RL}_H: \co^H \rightarrow \co_{H^{op}}.\]

\end{proposition}
 \begin{proof} The construction of $D^{RL}_H$ is similar to the construction of $D^{LR}_H$. The proof is an exercise.
\end{proof}



\begin{example} The dual of the standard comodule $(V, \D_V)$ on $\SLZ$.   \end{example}

Let us  define the \emph{left standard comodule} $(V, \D^l_V)$ on $\SLZ$ mentioned in Definition \ref{standard} as the following  map:

\[ \D^l_V: V \rightarrow  V \otimes R .\]

Let $V:=\Z\{e_1,e_2\}$ and $V^*:=\Z\{x_1,x_2\}$ with $x_i:=e_i^*$ the dual coordinates, and define 

\[ \D^l_{V}(e_1):= e_1\otimes x_{11}   + e_2\otimes x_{21} , \D^l_{V}(e_2):= e_1\otimes x_{12}   + e_2 \otimes x_{22} .\]

One checks we get the following: The map $\D_{V,1}$ is the following map

\[ \D_{V,1}(x_1):= x_{11} \otimes x_1 + x_{12}\otimes x_2, \D_{V,1}(x_2):= x_{21} \otimes x_1 + x_{22} \otimes x_2 .\]

The involution $S$ gives us the following dual comodule structure:

\[ \D_{V^*}: V^* \rightarrow R^{op} \otimes V^* \]

defined by 

\[ \D_{V^*}(x_1):= x_{22} \otimes x_1 - x_{12}\otimes x_2, \D_{V^*}(x_2):= -x_{21} \otimes x_1 + x_{11} \otimes x_2 .\]

We get the same comodule as in Example \ref{dualgroup}.

We may form the tensor product comodule

\[ \D_2: V^*\otimes_k V \rightarrow V^*\otimes_k V \otimes_k R \]

and there  are trace maps $tr: V^*\otimes V \rightarrow \Z$ and $tr \otimes 1: V^*\otimes V \otimes R \rightarrow R$ and one checks there is a commutative diagram. Hence with this definition it follows
the "pairing" $tr: V^* \otimes V \rightarrow \Z$ is invariant under the coaction of the Hopf algebra. Hence the functor $D^{LR}(-)$ gives a definition of the dual comodule for comodules that are free and of finite rank over $\Z$.

Note: Any torsion free $\Z$-module is free hence we may define the dual comodule for any  torsion free comodule. Note also that the "classical dual" from Example \ref{naivedual} differs from the dual from Definition \ref{dual} 
in the sense that the dual is a comodule on the opposite group scheme $\SLZ^{op}$. The "naive dual" is a comodule on $\SLZ$.

\begin{proposition} Let $(R, \D_R, S, \e)$ be a Hopf algebra with $H:=\Spec(R)$. Let $(V, \D) \in \co^H$ be a comodule on $R$ with $V$ a finite rank free $k$-module. Let $(V^{**}, \D^{**}):=(D^{RL}_H\circ D^{LR}_{H^{op}})(V, \D)$.
It follows $V^{**}$ is the double dual of $V$ and the canonical isomorphism $\phi: V \cong V^{**}$ gives an isomorphism of comodules 

\[  (V, \D) \cong (V^{**}, \D^{**}).\]

 Similarly if $(W,\D_1) \in \co_H$ is a comodule and if 
$(W^{**}, \D_1^{**}):=(D^{LR}_H\circ D^{RL}_{H^{op}})(W, \D_1)$ it follows there is a canonical isomorphism 

\[  (W, \D_1) \cong (W^{**}, \D_1^{**}) \]

 of comodules. Here $W^{**}$ is the double dual of $W$.
\end{proposition} 
\begin{proof} Let $\D_1: V \rightarrow R \otimes V$ be a comodule with $V:=k\{e_1,..,e_n\}$ and $V^*:=k\{x_1,..,x_n\}$ with $x_i:=e_i^*$. Let $x_i^*$ be the dual basis of  the basis $\{x_i\}$. Assume

\[ \D_1(e_i):= a_{i1} \otimes e_1+ a_{i2} \otimes e_2 + \cdots + a_{in} \otimes e_n.\]

By construction  it follows the map $\D_1^*: V^* \rightarrow V^* \otimes R^{op}$ is defined as follows:

\[ \D_1^*(x_i):= x_1 \otimes S(a_{1i}) + x_2\otimes S(a_{2i}) + \cdots + x_x \otimes S(a_{ni}).\]

It follows we get a comodule $\D_1^{**}: V^* \rightarrow R \otimes V^{**}$ defined by

\[ \D_1^{**}(x_i^*):= S^2(a_{i1}) \otimes x_1^* + S^2(a_{i2})\otimes x_2^* + \cdots + S^2(a_{in})\otimes x_n^* =\]

\[  a_{i1} \otimes x_1^* + a_{i2}\otimes x_2^* + \cdots +  a_{in}\otimes x_n^*. \]

It follows the canonical map $\phi: V \rightarrow V^{**}$ is an isomorphism of comodules, and the first claim of the Proposition follows. The second claim is similar and the Proposition is proved.

\end{proof}

\begin{definition}\label {dual}  Let $(V, \D_V)\in \co_H$ be a comodule on $R$ and let 

\[  (V^*, \D_{V^*}):= D^{LR}_H(V, \D_V)  \]

 be the comodule on $R^{op}$ from Proposition \ref{functor2}. We let the pair $(V^*, \D_{V^*})$ be the \emph{right dual comodule} of $(V, \D_V)$. Let $(W, \D_W) \in \co^H$ and define 

\[ (W^*, \D_{W^*}):= D^{RL}_H(W ,\D_W).   \]

We let the pair $(W^*, \D_{W^*})$ be the \emph{left dual comodule} of $(W, \D_W)$.
\end{definition}

Note: The dual $(V, \D_{V^*})$ is a comodule on the opposite Hopf algebra $R^{op}$ (or opposite group scheme $H^{op}$).

When we compose the functors $D^{LR}_H$ and $D^{RL}_H$ we get canonical isomorphisms of comodules, hence it makes sense to use these functors to define the dual comodule.  The definition may be used to define the dual for any comodule $(V, \D_V) $ on $\SLZ$ where $V$ is a finitely generated and torsion free $\Z$-module.

\begin{example} The dual of the symmetric square. \end{example}

Let $\SLZ:=\Spec(R)$ and let $\D: V \rightarrow R \otimes V$ be the following comodule: Let $V:=\Z\{e_1,e_2\}$ with $B:=\{e_1,e_2\}, C:=\{1\otimes e_1, 1\otimes e_2\}$ as bases and with the following comodule map:

\[ \D(e_1):= x_{11}\otimes e_1 + x_{12}\otimes e_2, \D(e_2):= x_{21}\otimes e_1 + x_{22} \otimes e_2.\]

Here we view $V$ as a free rank 2 module on $\Z$ and an element $u \in V$ expressed in the basis $B$  is a "vector" on the following form: $u:=u_1e_1+u_2e_2$ has $[u]_B:=[u_1,u_2]$.
Let 

\[
 M=  
\begin{pmatrix}   x_{11} & x_{12}  \\
x_{21}  &   x_{22}
\end{pmatrix}
\]

It follows 

\[  [\D ]^B_C([u]_B) = [\D(u)]_C= [u]_BM .\]

Hence we multiply the vector $[u]_B:=[u_1,u_2]$ from the right with the matrix $M$.

We may construct the symmetric power $\D_2: \Sym^2(V) \rightarrow \Sym^2(V) \otimes R^{op}$ as follows. Let $\Sym^2(V) \cong \Z\{e_1^2, e_1e_2, e_2^2\}$ with $B_2:=\{e_1^2, e_1e_2, e_2^2\}$ and let 
$C_2:=\{  e_1^2\otimes 1,  e_1e_2 \otimes 1,  e_2^2\otimes 1\}$.   Define the following matrix:

\[
 M_2=  
\begin{pmatrix}   x_{11}^2 & x_{11}x_{21} & x_{21}^2   \\
2x_{11}x_{12}   &   x_{11}x_{22}+x_{12}x_{21} & 2x_{21}x_{22} \\
x_{12}^2 & x_{12}x_{22}        &  x_{22}^2
\end{pmatrix}.
\]

Write elements $v \in \Sym^2(V)$ as column vectors. It follows 

\[  [\D_2]^{B_2}_{C_2}[v]_{B_2} =[\D_2(u)]_{C_2}=  M_2[v]_{B_2}.\] Hence we multiply the vector $[v]_{B_2}$ from the left with the matrix $M_2$.

Using the methods introduced in this section, we may construct the dual of the symmetric power 

\[  \D_3: \Sym^2(V)^* \rightarrow \Sym^2(V)^* \otimes R^{op} \]

 as follows. Let $\Sym^2(V)^* \cong \Z\{(e_1^2)^*, (e_1e_2)^*, (e_2^2)^*\}$ with $B_3:=\{(e_1^2)^*, (e_1e_2)^*,(e_2^2)^*\}$ and let 
$C_3:=\{ ( e_1^2)^*\otimes 1,  (e_1e_2)^* \otimes 1, ( e_2^2)^*\otimes 1\}$.   Define the following matrix:

\[
 M_3=  
\begin{pmatrix}   x_{22}^2 & -2x_{12}x_{22} & x_{12}^2   \\
-x_{21}x_{22}   &   x_{11}x_{22}+x_{12}x_{21} & -x_{11}x_{12} \\
x_{21}^2 & -2x_{11}x_{21}        &  x_{11}^2
\end{pmatrix}.
\]

Write elements $w \in \Sym^2(V)^*$ as column vectors. It follows 

\[  [\D_3]^{B_3}_{C_3}[w]_{B_3} =[\D_3(w)]_{C_3}=  M_3[w]_{B_3}.\] Hence we multiply the vector $[w]_{B_3}$ from the left with the matrix $M_3$.

You may check explicitly that the weight space decomposition of $\Sym^2(V)$ and $\Sym^2(V)^*$ agree, hence we cannot distinguish these two representations using weights. 
You may also check explicitly there is no isomorphism of comodules $\phi: \Sym^2(V) \cong \Sym^2(V)^*$. Hence the  two comodules $ \Sym^2(V)$ and $\Sym^2(V)^*$ are not isomorphic as comodules 
on $\SLZ^{op}$. Similar results
hold for higher symmetric powers: The comodules  $\Sym^d(V)$ and $\Sym^d(V)^*$ have the same weights and ranks of weight spaces,  but they are not isomorphic as comodules in general.

Let $B:=A[x_{ij}](D-1)$ with $\Spec(B):=\SLA$ and let $V:=A\{e_1,e_2\}$ be the following comodule on $\SLA$: Let $\D: V \rightarrow R\otimes_A V$ be defined by

\[ \D(e_1):=x_{11}\otimes e_1 +x_{12} \otimes e_2, \D(e_2):= x_{21} \otimes e_1 + x_{22} \otimes e_2.\]

\begin{proposition} Let $\Sym^2(V)$ be the symmetric square of $V$ and let $\Sym^2(V)^*$ be its dual. There is an isomorphism of comodules $\Sym^2(V) \cong \Sym^2(V)^*$ iff $2$ is a unit in $B$. Hence the comodules are isomorphic
when $A \cong \Z[\frac{1}{2}]$.
\end{proposition}
\begin{proof} The proof follows from the above construction and an explicit calculation.
\end{proof}






\section{The second symmetric tensor of the standard comodule}

In this section we prove that the second symmetric tensor $\Sym_2(V)$ and its dual $\Sym_2(V)^*$ have the same weights and weight spaces but are non-isomorphic as comodules on $\SLZ$. We also prove they are isomorphic as comodules
on $\SL_{2,\Z[\frac{1}{2}]}$. Hence a torsion free finite rank comodule on $\SLZ$ is not uniquely determined by its weights and weight space decomposition.

Let $R:=\Z[x_{ij}]/(D-1)$ be the Hopf algebra with the property that for any commutative $\Z$-algebra $B$ it follows the functor of points $h_R(B)$  equals the group of 2 by 2 matrices with coefficients in$B$ and determinant one.
It follows $G:=\Spec(R)$ equals the group scheme $G\cong \SLZ$ studied earlier in the paper.
Let $V:=\Z\{e_1,e_1\}$ be the free rank two $\Z$-module on the  elements $e_1,e_2$ and let $C:=\{e_1,e_2\}$. Let $V^*:=\Z\{x_1,x_2\}$ with $x_i:=e_i$ the dual basis and let $C^*:=\{x_1,x_2\}$.
Write $V\otimes_{\Z}R:=\{ e_1 \otimes 1, e_2 \otimes 1\}R$ as right $R$-module with basis $C_1:=\{e_1 \otimes 1, e_2 \otimes 1\}$.

Let $\D: V \rightarrow V \otimes_{\Z} R$ be the following comodule: Define

\[ \D(e_1):= e_1 \otimes x_{11}+ e_2 \otimes x_{21}, \D(e_2):= e_1 \otimes x_{12} + e_2 \otimes x_{22}.\]

\begin{example} The $d$'th symmetric tensors $\Sym_d(V)$. \end{example}

Let $d \geq 2$ be an integer. The symmetric group $S_d$ on $d$ letters acts on the $d$'th tensor power $V^{\otimes_{\Z} d}$ and we define the \emph{$d$'th symmetric tensors} $\Sym_d(V) \subseteq  V^{\otimes_{\Z} d}$ as the invariants
under this action. There is a comodule structure

\[  \D_d: \Sym_d(V) \rightarrow \Sym_d(V) \otimes R \]

and one may check that $\Sym_d(V)$ and $\Sym^d(V)$ have the same weights and the same ranks of weight spaces.


Note: When studying group schemes and representations of group schemes one prefer to study the symmetric tensors $\Sym_d(V)$ over the symmetric powers $\Sym^d(V)$, since the symmetric  tensors have better
functorial properties. There is over a field $k$ of characteristic zero an isomorphism of representations $\Sym_d(V\otimes k) \cong \Sym^d(V\otimes k)$. There is an equality of classes $[\Sym_d(V)]= [\Sym^d(V)]$ in $\operatorname{K}(\SLZ)$, hence the grothendieck group $\operatorname{K}(\SLZ)$ is generated by the classes $[\Sym_d(V)]$ of the symmetric tensors for $d\geq 0$ an integer.

In this section we will study the symmetric tensor $\Sym_2(V)$ and its dual $\Sym_2(V)^*$. 


Let

\[
 M=  
\begin{pmatrix}   x_{11} & x_{12}  \\
         x_{21}   &   x_{22}
\end{pmatrix}.
\]

Write the vector $u:=u_1e_1+u_2e_2$ in the basis $C$ as a column vector

\[
 [u]_C=  
\begin{pmatrix}   u_1 \\
          u_2  
\end{pmatrix}.
\]

It follows the map $\D$ is the following map:

\[  [\D(u)]_{C_1} = M[u]_C .\]

Hence if we view the map $\D$ as a map from $V$ in the basis $C$ to $V \otimes R$ in the basis $C_1$ it follows $\D$ is given as left multiplication with the matrix $M$. Hence there is an equality

\[  M= [\D ]^C_{C_1} .\]

There is an  isomorphism $\Sym_2(V) \cong \Z\{ e_{11}, e_{12}, e_{22} \}$ where $e_{11}:=e_1 \otimes e_1, e_{12}:=e_1 \otimes e_2 + e_2 \otimes e_1$ and  $e_{22}:=e_2 \otimes e_2$. Let $E:=\{e_{11}, e_{12}, e_{22}\}$. 
Let $\Sym_2(V) \otimes R \cong \{e_{11} \otimes 1, e_{12} \otimes 1, e_{22} \otimes 1\}R$ and let $E_1:=\{e_{11} \otimes 1, e_{12} \otimes 1, e_{22} \otimes 1\}$.

We get a comodule structure $\D_2: \Sym_2(V) \rightarrow \Sym_2(V) \otimes R$ defined as follows:

\[ \D_2(e_{11}):= e_{11} \otimes x_{11}^2 + e_{12} \otimes x_{11}x_{21} + e_{22} \otimes x_{21}^2, \]

\[ \D_2(e_{12}):= e_{11} \otimes 2x_{11}x_{12} + e_{12} \otimes (x_{11}x_{22}+ x_{12}x_{21}) + e_{22} \otimes 2x_{21}x_{22}, \]

and

\[ \D_2(e_{22}):= e_{11} \otimes x_{12}^2 + e_{12} \otimes x_{12}x_{22}  + e_{22} \otimes x_{22}^2. \]

Let a vector $u:=u_1e_{11}+ u_2 e_{12} + u_3 e_{22}$ be written as a column vector in the basis $E$:

\[
 [u]_E=  
\begin{pmatrix}   u_2 \\
          u_2 \\
u_3 
\end{pmatrix}.
\]

We may write the map $\D_2$ in the bases $E, E_1$. Let $M$ be the following matrix:

\[
 M:=  
\begin{pmatrix}   x_{11}^2 & 2x_{11}x_{12} & x_{12}^2  \\
            x_{11}x_{21} & x_{11}x_{22} + x_{12}x_{21} & x_{12}x_{22} \\
 x_{21}^2 & 2x_{21}x_{22} & x_{22}^2 
\end{pmatrix}.
\]

It follows there is an equality

\[  [\D_2(u)]_{E_1}= M[u]_E.\]

Let $\Sym_2(V)^* \cong \Z\{y_{11}, y_{12}, y_{22} \}$ with $y_{ij}:=e_{ij}^*$ the dual basis, and let $E^*:=\{ y_{ij} \}$. Represent a vector $v:= v_1 y_{11}+ v_2y_{12} + v_3y{22}$ in the basis $E^*$ as a row vector 

\[ [v]_{E^*}:=[v_1,v_2,v_3] .\]

Let $ R^{op} \otimes \Sym_2(V)^* \cong R^{op}\{ 1 \otimes y_{ij} \}$ and let $E_1^*:=\{ 1\otimes y_{ij} \}$.

Define the following map:

\[ \D_2^*: \Sym_2(V)^* \rightarrow R^{op} \otimes \Sym_2(V)^* \]

by

\[ \D_2^*(v):=  S\otimes 1 \circ \tau_2 \circ \D_2(v) \]

where

\[  \tau_2: \Sym_2(V)^* \otimes R \rightarrow R \otimes \Sym_2(V)^*  \]

is the map defined by $\tau_2(v \otimes z):= z \otimes v $ from Lemma \ref{involution}. 
Define the map $S\otimes 1: R\otimes \Sym_2(V)^* \rightarrow R \otimes \Sym_2(V)^*$ using the involution $S$ of $R$. From Lemma \ref{involution} and Proposition \ref{functor2} it follows $(\Sym_2(V)^*, \D_2^*)$ is a comodule on $R^{op}$. 
Let $M^*$ be the following matrix:

\[
 M^*:=  
\begin{pmatrix}   x_{22}^2 & -2x_{12}x_{22} & x_{12}^2  \\
            -x_{21}x_{22} & x_{11}x_{22} + x_{12}x_{21} & -x_{11}x_{12} \\
 x_{21}^2 & -2x_{11}x_{21} & x_{11}^2 
\end{pmatrix}.
\]

It follows we may describe the map $\D_2^*$ is the bases $E^*, E^*_1$ using the matrix $M^*$: There is an equality

\[ [\D_2^*(v)]_{E^*_1} = [v]_{E^*} M^*. \]

Hence we multiply the vector $[v]_{E^*}$ from the right with the matrix $M^*$. Hence the map $[\D_2^*]^{E^*}_{E_1^*}$ is given as follows: We multiply a vector with the matrix $M^*$ from the right.

We may using Lemma \ref{involution} apply the involution $\tau_2$ to the comodule $\D_2: \Sym_2(V) \rightarrow \Sym_2(V) \otimes R$ to get a comodule on the opposite Hopf algebra $R^{op}$. Let $\D_2^{\tau_2}:= \tau_2 \circ \D_2$. We get a comodule

\[ \D_2^{\tau_2}: \Sym_2(V) \rightarrow R^{op} \otimes \Sym_2(V).\]

If  we write a vector $u:=u_1e_{11} + u_2 e_{12} + u_3 e_{22}$  in the basis $E$, let $[u]_E^{tr}:=[u_1,u_2,u_3]$ denote the transpose vector, which is a row vector. When we use the involution $\tau_2$ we take the transpose of matrices.
Let $R^{op}\otimes \Sym_2(V) \cong R^{op}\{1 \otimes e_{ij} \}$ and let $E_2:=\{1\otimes e_{ij} \}$. In the bases $E,E_2$ we may write the map $\D_2^{\tau_2}$ as follows:

\[ [\D_2^{\tau_2}(u)]_{E_2}  = [u]_E^{tr} M^{tr}.\]

Here $M^{tr}$ is the transpose of the matrix $M$ and we multiply the row vector $[u]_E^{tr}$ from the right with the matrix $M^{tr}$.
The matrix $M^{tr}$ is the following matrix:

\[
 M^{tr}:=  
\begin{pmatrix}   x_{11}^2 & x_{11}x_{21} & x_{21}^2  \\
           2 x_{11}x_{12} & x_{11}x_{22} + x_{12}x_{21} & 2x_{21}x_{22} \\
 x_{12}^2 & x_{12}x_{22} & x_{22}^2 
\end{pmatrix}.
\]

Let $A \in \GL_{\Z}(\Z^3)$ be a matrix with
\[
 A:=  
\begin{pmatrix}   a_{11} & a_{12}  & a_{13}  \\
           a_{21} & a_{22}  & a_{23}    \\
 a_{31} & a_{32}  & a_{33}
\end{pmatrix}.
\]

One checks there is an isomorphism of comodules  $\phi: \Sym_2(V) \cong \Sym_2(V)^*$ iff there is a matrix $A\in \GL_{\Z}(\Z^3)$ with  $AM^{tr}=M^* A$. This gives a system of equations in the "variables" $a_{ij}$ and one finds $\phi$ exists iff the matrix $A$  has the following form:

\[
 A:=  
\begin{pmatrix}   0  & 0   & -2a  \\
           0  & a   &   0    \\
 -2a &  0   &  0
\end{pmatrix}.
\]

If $a=1$ it follows $det(A)=4$ hence $A$ is invertible iff $2$ is a unit. Hence such a matrix exsists in $\GL_{\Z[\frac{1}{2}] }(\Z[\frac{1}{2}])^3)$. There is no such matrix in $\GL_{\Z}(\Z^3)$.

\begin{theorem} The comodules $\Sym_2(V)$ and $\Sym_2(V)^*$ have the same weights and ranks of weight spaces but are not isomorphic as comodules on $\SLZ$. The pull backs $\Sym_2(V)\otimes_{\Z} \Z[\frac{1}{2]}]$ and 
$\Sym_2(V)^* \otimes_{\Z} \Z[\frac{1}{2}]$ are isomorphic as comodules on $\SL_{2,\Z[\frac{1}{2}]}$.
\end{theorem}
\begin{proof} One checks the comodules have the same weights and ranks of weight spaces and an explicit calcualtion proves they are not isomorphic as comodules on $\SLZ$ using the above constructions. One also proves they are isomorphic
as comodules on $\SL_{2,\Z[\frac{1}{2}]}$ using an explicit calculation with the matrices $M^*$ and $M^{tr}$. The Theorem follows.
\end{proof}

\begin{example}  \label{nogoodfiltration} Comodules on $\SLZ$ and $\SL_{2, \mathbb{Q}}$. \end{example}

The comodules $\Sym^2(V)^*$ and $\Sym_2(V)^*$ do not have a good filtration. Hence it is easy to give explicit examples of comodules on $\SLZ$ not having a good filtration.

Let $S:=\Spec(\Z)$ and $T:=\Spec(k)$ with $k$ the field of rational numbers. We may view the group scheme $G_k:=\SL(2,k)$ as the special fiber
of the canonical map $\pi: \SLZ \rightarrow S$ and we may pull back any comodule on $\SLZ$ to $G_k$. Let $\tau: G_k \rightarrow \SLZ$ be the canonical map.
Given any torsion free finite rank comodule $W$ on $\SLZ$ we get a canonical decomposition 

\[ \tau^*(W) \cong \oplus_{i=1}^l \Sym^{d_i}(k^2) \]

where $k^2$ is the standard comodule on $\SL(2,k)$. We would like for each integer $d \geq 0$ a classification of the comodules $W_d$ on $\SLZ$ with $\tau^*(W_d) \cong \Sym^d(k^2)$.
The methods introduced above gives a method to classify all such comodules. As an example: 

\begin{theorem} \label{descent} The set of finite rank torsion free comodules on $\SLZ$ with $\tau^*(W) \cong \Sym^2(k^2)$ are the following modules:

\[  \Sym^2(V), \Sym^2(V)^*, \Sym_2(V), \Sym_2(V)^*.\]

\end{theorem}
\begin{proof} This is an explici calculation using the methods introduced in this section. 
\end{proof}

\section{A universal Clebsch-Gordan filtration for $\GL_{2,A}$ }

In this section I prove the existence of a finite filtration $E_i$ of $\GL_{2,A}$-comodules of the tensor product $\Sym^n_A(V) \otimes \Sym^m_A(V)$ for all integers $1 \leq n \leq m$ and. Here $V$ is the standard comodule.
The filtration $E_i$ i a good filtration, hence the construction gives an explicit and elementary construction of an infinite class of finite rank torsion free comodules on $\GL_{2,A}$ with a good filtration.
When we restrict the filtration $E_i$to the closed subgroup $\SL_{2,A}$ we get the filtration $F_i$ defined earlier.

Note: The proof is done in terms of the functor of points $h_G(-)$ of $G$ and not in terms of the comodule $k[G]$.

Let in this section $A$ be any commutative ring and let $V:=A\{e_1,e_2\}$ be the free $A$-module of rank two on the elements $e_1,e_2$. Let $B:=\{e_1,e_2\}$ and for any $R\in \alg$ let $B(R):=\{1_R \otimes e_1, 1_R \otimes e_2\}$
be the basis for $R\otimes_A V$. Let $e(i)_R:=1_R \otimes e_i$. 
Let $G:=\GL_{2,A}$ have coordinate ring $k[G]:=A[x_{ij}, 1/D]$
where $x_{11}, x_{12}, x_{21}, x_{22} $ are independent variables over $A$ and $D:=x_{11} x_{22}- x_{12} x_{21} $ is the determinant.  Let $C:=\{x_{11}, x_{12}, x_{21}, x_{22} \}$ and for any element $g \in h_G(R)$ with $R \in \alg$ let 

\[
 [g]_C:=  
\begin{pmatrix}   g(x_{11})  & g(x_{12})  \\
            g(x_{21})  & g(x_{22})        \\
\end{pmatrix}.
\]

Let $h_V$ be the following functor:

\[ h_V(-): \alg \rightarrow \grp \]

defined by 

\[ h_V(R):=\{R\otimes_A V, +\} \]

where $\{R\otimes_A V,+ \}$ is the underlying additive abelian group of the $A$-module $R\otimes_A V$. Obviously $h_V(-)$ is a well defined functor.  Define for any $R\in \alg$ the following action:

\[  <,>^1_R: h_G(R) \times h_V(R) \rightarrow h_V(R) \]

by

\[  <g,   u>^1_R:=   [g]_C [u]_{B(R)} =\]

\[
 \begin{pmatrix}   g(x_{11})  & g(x_{12})  \\
            g(x_{21})  & g(x_{22})        \\
\end{pmatrix}
\begin{pmatrix}    u_1 \\
                            u_2 \\
\end{pmatrix}=
\begin{pmatrix}   g(x_{11})u_1  + g(x_{12})u_2  \\
            g(x_{21})u_1+ g(x_{22}) u_2       \\
\end{pmatrix}.
\]

where $u:=u_1e(1)_R + u_2e(2)_R$.

\begin{lemma} \label{leftaction} The following holds: For any $a\in R, u,v \in u_V(R)$ it follows

\[ <g,u+v>^1_R=<g,u>^1_r+<g,v>^1_R, <g,au>^1_R=a<g,u>^1_R .\]

Let $g,h\in h_G(R)$ and $e\in h_G(R)$ be the identity. It follows

\[  <g, <h,u>^1_R>^1_R=<gh, u>^1_R,  <e,u>^1_R=u. \]

\end{lemma}
\begin{proof} The proof follows directly from the definition. \end{proof}

\begin{corollary} For any $R\in \alg$ we get a map of groups

\[ \rho_R: h_G(R) \rightarrow \GL_R(R\otimes_A V) \]

defined by $\rho_R(g)(u):=<g,u>^1_R$. For any map $f:R \rightarrow R_1 $ in $\alg$  we get a commutative diagram

\[
\diagram     h_G(R) \rto \dto &  h_V(R) \dto  \\
                   h_G(R_1) \rto & h_V(R_1) 
\enddiagram.
\]
Hence we get a natrual tranformation $\rho: h_G(-) \rightarrow h_V(-)$ of functors and a representation of $\GL_{2,A}$ in $\GL_A(V)$.
\end{corollary}

\begin{definition} Let $\rho: \GL_{2,A} \rightarrow \GL_A(V)$ be the \emph{standard representation} of $\GL_{2,A}$.
\end{definition}

The aim of this section is to study the standard representation $(V,\rho)$ and various Schur modules of $V$. 

Let $f:R\rightarrow R_1 \in \alg$ be a map and consider the symmetric power $\Sym^d_A(V)$ and exterior power $\wedge^d_A V$  for an integer $d \geq 1$. There is a canonical map

\[ f\otimes 1: R\otimes_A \Sym^d_A(V) \rightarrow R_1\otimes_A \Sym^d_A(V)  \]

and similar for $R\otimes_A \wedge^d_AV \rightarrow R_1 \otimes_A \wedge^d_A V$. There  are canonical induced maps

\[  h_{\Sym^d(V)}(f): \Sym^d_R(R\otimes_A V) \rightarrow \Sym^d_{R_1}(R_1 \otimes_A V)  \]

and

\[  h_{\wedge^d V}(f): \wedge^d_R(R\otimes_A V ) \rightarrow \wedge^d_{R_1}(R_1 \otimes_A V) \]

defined by 

\[ h_{\Sym^d(V)}(f) (x_1 \otimes v_1 ) \cdots (x_d \otimes v_d):= (f(x_1)\otimes v_1) \cdots (f(x_d) \otimes v_d) \]

and 

\[ h_{\wedge^d V}(f)(x_1 \otimes v_1)\wedge \cdots \wedge (x_d \otimes v_d):=(f(x_1) \otimes v_1) \wedge \cdots \wedge (f(x_d) \otimes v_d).\]

\begin{lemma} Let $f:R\rightarrow R_1 \in \alg$ be a map and consider $\Sym^d_A(V)$ and $\wedge^d_A V$  for an integer $d \geq 1$. There are canonical isomorphisms

\[  \phi_R: R\otimes_A \Sym^d_A(V) \cong \Sym^d_R(R\otimes_A V), \psi_R: R\otimes_A \wedge^d_A V \cong \wedge^d_R(R \otimes_A V), \]

where 

\[ \phi_R( x_1 \otimes v_1v_2 \cdots v_d):=(x_1\otimes v_1)(1_R \otimes v_2) \cdots (1_R \otimes v_d).\]

Similarly for $\psi_R$. There is a commutative diagram 

\[
\diagram     R\otimes_A \Sym_A^d(V) \rto^{\phi_R}    \dto^{f\otimes 1}   &  \Sym^d_R(R\otimes_A V) \dto^{h_{\Sym^d(V)}(f) }  \\
                   R_1 \otimes_A \Sym^d_A(V) \rto^{\phi_{R_1} } & \Sym^d_{R_1}(R_1\otimes_A V) 
\enddiagram.
\]

There is a similar comutative diagram for the exterior power.
\end{lemma} 
\begin{proof} The proof is immediate.
\end{proof}

Let $g\in h_G(R)$ and consider the following map

\[  <g, ->: (R\otimes_A V)^{\times d}  \rightarrow \Sym^d_R(R\otimes_A V)  \]

defined by 

\[ <g, (u_1,u_2,..,u_d)>:= (<g,u_1>^1_R)(<g,u_2>^1_R) \cdots (<g,u_d>^1_R) .\]

It follows the map $<g,->$ i $R$-multilinear and hence it induce a map

\[ <g,->^d_R: \Sym^d_R(R\otimes_A V) \rightarrow \Sym^d_R(R\otimes_A V) \]

defined by 

\[  <g, \sum_i u(i)_1\cdots u(i)_d>^d_R:=\sum_i (<g,u(i)_1>^1_R)(<g,u(i)_2>^1_R) \cdots (<g,u(i)_d>^1_R) .\]

There is a commutative diagram

\[
\diagram     R\otimes_A \Sym_A^d(V) \rto^{\phi_R}    \dto^{f\otimes 1}   &  \Sym^d_R(R\otimes_A V) \dto^{h_{\Sym^d(V)}(f) }  \\
                   R_1 \otimes_A \Sym^d_A(V) \rto^{\phi_{R_1} } & \Sym^d_{R_1}(R_1\otimes_A V) 
\enddiagram.
\]


Let  $f: R \rightarrow R_1 \in \alg$ be any map.  Let $G(f): h_G(R) \rightarrow h_G(R_1)$ be the induced map. Let $h_{\Sym^d(V)}(f):=h^d(f)$. The  action 

\[ <,>^d_R: h_G(R) \times \Sym^d_R(R\otimes_A V) \rightarrow \Sym^d_R(R\otimes_A V) \]

has the following properties: Let $a\in R, u,v \in R\otimes_A V, g,h,e \in h_G(E)$. It follows

\[ <g,u+v>^d_R=<g,u>^d_R + <g,v>^d_R, <g,au>^d_R =a<g,u>^d_R \]

and

\[  <g,<h,u>^d_R>^d_R=<gh, u>^d_R, <e,u>^d_R=u.\]

 The above claims hold because of properties of the bracket $<-,->^1_R$.

There is a commutative diagram

\[
\diagram     G(R) \times  \Sym^d_R(R\otimes_A V) \dto^{G(f) \times h^d(f)  }  \rto &  \Sym^d_R(R\otimes_A V) \dto^{h^d(f)}   \\
                   G(R_1) \times  \Sym^d_{R_1}(R_1\otimes_A V) \rto  & \Sym^d_{R_1}(R_1 \otimes_A V) 
\enddiagram.
\]

This claim follows from the formula

\[ <G(f)(g), h^1(f)(u)>^1_{R_1} =h^1(f)(<g,u>^1_R)  \]

which holds for all $f,g,u$.

For every integer $d\geq 1$ and $R$  we get a map of groups

\[ \rho^d_R: G(R) \rightarrow \GL_R(\Sym^d_R(R\otimes_A V)) \]

defined by $\rho^d(g)(u):=<g,u>^d_R$ and a map of functors 

\[ \rho^d: \GL_{2,A} \rightarrow \GL_A(\Sym^d_A(V)).  \] 

Hence we may view $\{\Sym^d_A(V), \rho^d\}$ as a representation of the group scheme $\GL_{2,A}$.


\begin{definition} The pair $\{\Sym^d_A(V), \rho^d\}$ is the \emph{$d$'th symmetric power of $V$}.
\end{definition}

We may do the same for the exterior power (or for any Scur-Weyl functor): For any $R \in \alg, g\in h_G(R)$ and integer $d\geq 1$ define an action

\[  <g, ->: \wedge^d_R(R\otimes_A V) \rightarrow \wedge^d_R(R\otimes_A V) \]

by

\[  <g, (x_1 \otimes v_1)\wedge \cdots \wedge (x_d \otimes v_d)>:=(<g, x_1 \otimes v_1>^1_R)\wedge \cdots  \wedge (<g,x_d\otimes v_d>^1_R) .\]

We get a representation

\[ \eta^d: \GL_{2,A} \rightarrow \GL_A(\wedge_A^d V).\]

\begin{definition} The pair $\{ \wedge^d_A V, \eta^d\}$ is the \emph{$d$'th exterior power of $V$}.
\end{definition}

For any pair of $\GL_{2,A}$-modules $W_1, W_2$ we may define the tensor product $W_1 \otimes_A W_2$  and the functor

\[ h_{W_1\otimes W_2}: \alg \rightarrow \grp \]

by

\[ h_{W_1 \otimes W_2}(R):=\{R\otimes_A W_1 \otimes_A W_2, +\}.\]

There is a canonical isomorphism $h_{W_1 \otimes_A W_2}(R) \cong (R\otimes_A W_1)\otimes_R (R\otimes_A W_2) $ and for any $R\in alg$ an action

\[ <-, ->: G(R) \times h_{W_1 \otimes W_2}(R) \rightarrow h_{W_1 \otimes W_2}(R) \]

defined by 

\[  <g, u_1 \otimes u_2>:= <g, u_1>_{R,1} \otimes < g, u_2>_{R,2} \]

where $<-,->_{R,i}$ is the action on $R\otimes_A W_i$. 

We get a representation 

\[  \gamma_{1,2}: \GL_{2,A} \rightarrow \GL_A(W_1\otimes_A W_2).\]


There is for every $R$ and integers $n,m \geq 1$ canonical map of $R$-modules

\[   \phi_{n,m}: \Sym^{n-1}_R(R\otimes_A V) \otimes \Sym^{m-1}_R(R\otimes_A V) \otimes \wedge^2_R(R\otimes_A V) \rightarrow   \Sym^n_R(R\otimes_A V) \otimes \Sym^m_R(R\otimes_A V) \]

defined by 

\[ \phi_{n,m}(z \otimes w \otimes u \wedge v ):=zu\otimes wv - zv\otimes wu , \]

and for any map $f: R \rightarrow R_1$ a commutative diagram

\[
\diagram      \Sym^{n-1}_R(R\otimes_A V) \otimes \Sym^{m-1}_R(R \otimes_A V)\otimes \wedge^2_R(R\otimes_A V) \rto \dto &  \Sym^n_R(R\otimes_A V) \otimes \Sym^m_R(R\otimes_A V) \dto   \\
   \Sym^{n-1}_{R_1}(R_1\otimes_A V) \otimes \Sym^{m-1}_{R_1}(R_1 \otimes_A V)\otimes \wedge^2_{R_1}(R_1\otimes_A V) \rto   &  \Sym^n_{R_1}(R_1\otimes_A V) \otimes \Sym^m_{R_1}(R_1\otimes_A V)              
\enddiagram.
\]

We get a map

\[ \psi_{n,m}: \Sym^n_R(R\otimes_A V) \otimes \Sym^m_R(R\otimes_A V) \rightarrow \Sym^{n+m}_R(R\otimes_A V) \]

defined by

\[ \psi_{n,m}(z \otimes w):=zw.\]

and a commutative diagram

\[
\diagram
\Sym^n_R(R\otimes_A V) \otimes \Sym^m_R(R\otimes_A V) \rto \dto & \Sym^{n+m}_R(R\otimes_A V)  \dto \\
\Sym^n_{R_1}(R_1\otimes_A V) \otimes \Sym^m_{R_1}(R_1\otimes_A V) \rto & \Sym^{n+m}_{R_1}(R_1\otimes_A V)
\enddiagram.
\]

\begin{lemma} \label{exact} For all integers $1 \leq n \leq m$ and all $R$ the following sequence is an exact sequence of left $R$-modules:

\[ 0 \rightarrow  \Sym^{n-1}_R(R\otimes_A V) \otimes \Sym^{m-1}_R(R\otimes_A V) \otimes \wedge^2_R(R\otimes_A V) \rightarrow \Sym^n_R(R\otimes_A V) \otimes \Sym^m_R(R\otimes_A V)  \rightarrow   \]

\[  \Sym^{n+m}_R(R\otimes_A V) \rightarrow 0.\]

\end{lemma}
\begin{proof} The proof is similar to the proof of Lemma \ref{zmodules}
\end{proof}

When $R:=A$ we get the following exact sequence of $A$-modules:

\[ 0 \rightarrow  \Sym^{n-1}_A(  V) \otimes \Sym^{m-1}_A(  V) \otimes \wedge^2_A V \rightarrow \Sym^n_A( V) \otimes \Sym^m_A( V)  \rightarrow    \Sym^{n+m}_A(  V) \rightarrow 0.\]

Using methods similar to the above methods we may prove the following result:

\begin{theorem} Let $1\leq n \leq m$ be integers. The sequence

\[ 0 \rightarrow  \Sym^{n-1}_A(  V) \otimes \Sym^{m-1}_A(  V) \otimes \wedge^2_A V \rightarrow \Sym^n_A( V) \otimes \Sym^m_A( V)  \rightarrow    \Sym^{n+m}_A(  V) \rightarrow 0\]

is an exact sequence of $\GL_{2,A}$-comodules.
\end{theorem}
\begin{proof}   Let $W_1:=\Sym^{n-1}_R(R\otimes_A V) \otimes \Sym^{m-1}_R(R\otimes_A V)\otimes \wedge^2_R(R\otimes_A V)$ and let $<-,->_1$ be the action of $h_G(R)$ on $W_R$. 
Let $W_2:=\Sym^n_R(R\otimes_A V) \otimes \Sym^m_R(R\otimes_A V)$ and let $<-,->_2$ be the action of $h_G(R)$ on $W_2$.

The map $\phi_{n,m}$ has the following property: $\phi_{n,m}(< g, z \otimes w \otimes u \wedge v>_1)=$

\[ \phi_{n,m}(g(z)\otimes g(w) \otimes g(u) \wedge g(v) )=\]

\[ g(z)g(u)\otimes g(w)g(v) - g(z)g(v)\otimes g(w)g(u)= g(zu\otimes wv - zv \otimes wu)=\]

\[ <g, \phi_{n,m}(z \otimes w \otimes u \wedge v) >_2\]

hence the map $\phi_{n,m}$ is a map of $h_G(R)$-modules. Similarly for $\psi_{n,m}$.

The Theorem follows.
\end{proof}

\begin{definition}

Let $E_i:=\Sym^{n-i}_A(V) \otimes \Sym_A^{m-i }(V) \otimes (\wedge^2_A V)^{\otimes_A^i}$. We get a filtration

\[0:= E_{n+1} \subsetneq E_n \subsetneq \cdots \subsetneq E_1 \subseteq E_0:=\Sym^n_A(V) \otimes \Sym^m_A(V) \]

called the \emph{universal Clebsch-Gordan filtration} for the $\GL_{2,A}$-comodule $\Sym_A^n(V)\otimes \Sym_A^m(V)$.

\end{definition}

There is for every $i$ a canonical isomorphism of comodules

\[ E_i/E_{i+1} \cong \Sym_A^{n+m-2i}(V) \otimes (\wedge^2_A V)^{\otimes_A i} .\]

In the Grothendieck group of finite rank $\GL_{2,A}$-comodules we get an equality of classes

\[   [\Sym^n_A(V)][\Sym_A^m(V)] =  \sum_{i=0}^n [\Sym_A^{n+m-2i}(V) ][\wedge^2_A V]^i ,\]

called the \emph{virtual Clebsch-Gordan formula} for $\GL_{2,A}$.

\begin{corollary}  If $A:=\Z$ and we restrict to $\SL_{2,\Z}$ we get the virtual Clebsch-Gordan formula for $\SL_{2,\Z}$ from Corollary \ref{virtualCG}. 
\end{corollary}
\begin{proof} The Corollary follows since $\wedge^2_{\Z}V$ is a trivial $\SL_{2,\Z}$-module.
\end{proof}

Note: In the previosu sections we gave a "constructive proof" of the virtual Clebsch-Gordan formula using comodules on Hopf algebras and the ring $k[\SL_{2,A}]$. The above proof is in some sense "non-constructive" since it uses the functor of points $h_{\GL_{2,A}}$,  "non small categories" and "universes".

\end{document}